\documentclass[12pt,onecolumn]{IEEEtran}
\usepackage{amssymb}
\usepackage{graphicx, epstopdf}
\usepackage{epstopdf}
\usepackage[tbtags]{amsmath}
\usepackage{amsfonts}
\usepackage{mathrsfs}
\usepackage{cite,url}
\usepackage{upquote}
\usepackage{pifont}
\usepackage{amsmath}
\usepackage{mathtools}

\usepackage{units}
\usepackage{xcolor}

\makeatletter
\newcommand{\vast}{\bBigg@{4}}
\newcommand{\Vast}{\bBigg@{5}}
\makeatother

\newtheorem{thm}{Theorem}
\newtheorem{theorem}[thm]{Theorem}
\newtheorem{lemma}[thm]{Lemma}
\newtheorem{definition}[thm]{Definition}
\newtheorem{proposition}[thm]{Proposition}
\newtheorem{corollary}[thm]{Corollary}

\newtheorem{example}{Example}

\newtheorem{remark}{Remark}

\begin{document}
	\title {Generalized Fisher-Darmois-Koopman-Pitman Theorem and Rao-Blackwell Type Estimators for Power-Law Distributions}
	\author{
		\IEEEauthorblockN{Atin Gayen and M. Ashok Kumar}\\
    \thanks{Atin Gayen is supported by an INSPIRE fellowship of the Department of Science and Technology, Govt. of India.}
    \thanks{Atin Gayen (Email: atinfordst@gmail.com) and M. Ashok Kumar (Email: ashokm@iitpkd.ac.in) are with the Department of Mathematics
of Indian Institute of Technology Palakkad, Kerala 678 557, India.}
\thanks{A part of this paper was presented at the 2021 IEEE International Symposium on Information Theory \cite{GayenK21ISIT}.}
	}
	
	\maketitle
	
\begin{abstract}
This paper generalizes the notion of sufficiency for estimation problems beyond maximum likelihood. In particular, we consider estimation problems based on Jones et al. and Basu et al. likelihood functions that are popular among distance-based robust inference methods. We first characterize the probability distributions that always have a fixed number of sufficient statistics (independent of sample size) with respect to these likelihood functions. These distributions
are power-law extensions of the usual exponential family and contain Student distributions as a special case. We then extend the notion of minimal sufficient statistics and compute it for these power-law families. Finally, we establish a Rao-Blackwell-type theorem for finding the best estimators for a power-law family. This helps us establish Cram\'er-Rao-type lower bounds for power-law families.
\end{abstract}

\begin{IEEEkeywords}
Cram\'er-Rao bounds,
Fisher-Darmois-Koopman-Pitman theorem,
Rao-Blackwell theorem,
Student distribution,
sufficient statistics.
\end{IEEEkeywords}

%
\IEEEpeerreviewmaketitle

\section{Introduction}
\label{2sec:introduction}
\IEEEPARstart{C}{onsider} a parametric family of probability distributions, parameterized by $\theta$, $\Pi = \{p_\theta:\theta\in\Theta\}$, where $\Theta$ is an open subset of $\mathbb{R}^k$. We shall assume (unless explicitly specified) that all distributions in the family have common support $\mathbb{S}\subset\mathbb{R}$. Let $X_1^n\coloneqq X_1,\dots, X_n$ be an independent and identically distributed (i.i.d.) random sample drawn according to some $p_\theta\in\Pi$. An estimator of $\theta$ is a function of $X_1^n$. If the estimator is a function of some specific statistics of the sample, then one can use these statistics to estimate $\theta$ without knowing the entire sample. Fisher introduced this concept in his seminal paper titled ``On the mathematical foundations of theoretical statistics'' and termed it as \emph{sufficient statistics} for the parameter \cite[Sec.~7]{Fisher22J}. He designated a statistic $T$ to be 
 \emph{sufficient for estimating $\theta$} if the conditional distribution of the sample, given $T$, is independent of $\theta$. Formally, if $T\coloneqq T(X_1^n)$ is a statistic ($T$ could very well be a vector-valued function), then $T$ is sufficient (for estimating $\theta$) if the conditional distribution
\begin{equation}
\label{2eqn:fisher_def}
    p_{\theta_{X_1^n|T}}(x_1^n|t)\quad\text{is independent of}\quad \theta
\end{equation}
for every $t\in \mathscr{T} \coloneqq\{T(y_1^n)\colon y_1^n\in \mathbb{S}^n\}$. This definition, however, is helpful only if one has a suitable guess for $T$. Observe also that the distribution of $T$ is needed to find the conditional distribution in (\ref{2eqn:fisher_def}). Later, Koopman proposed the following equivalent definition, where the distribution of $T$ is not needed \cite{Koopman36J}. He called a statistic $T$ sufficient for $\theta$ if the following holds. If $x_1^n$ and $y_1^n$ are two sets of i.i.d. samples such that they are same under $T$, that is, $T(x_1^n)=T(y_1^n)$, then the difference of their respective log-likelihoods, $[L(x_1^n;\theta) - L(y_1^n;\theta)]$ is independent of $\theta$, where $L(x_1^n;\theta) \coloneqq \sum_{i=1}^n \log p_\theta(x_i)$ denotes the log-likelihood function of the sample $x_1^n$. Observe that, to use this definition also, one needs to have a guess for $T$. Let us now suppose that a sufficient statistic $T$ is known. A natural question that arises is, what kind of estimators that can be found using $T$. Observe that both Fisher's and Koopman's definitions of sufficiency make use of the log-likelihood function $L$. However, it is still not clear if Fisher or Koopman had maximum likelihood (ML)-estimation in mind while proposing the above definitions. The following result, known as \emph{factorization theorem}  \cite{Neyman35J, HalmosS49J}, which is equivalent to both the above definitions, seems to answer this. A statistic $T$ is sufficient for $\theta$ if and only if there exist two functions $u:\Theta\times \mathscr{T}\to \mathbb{R}$ and $v:\mathbb{S}^n\to \mathbb{R}$ such that 
\begin{eqnarray}
    \label{2eqn:factorization_thm}
    L(x_1^n;\theta) = u(\theta,T(x_1^n)) + v(x_1^n)
\end{eqnarray}
for all $x_1^n$ and all $\theta\in\Theta$ \cite[Th. 6.2.6]{CasellaB02B}. Observe that, using this criterion, one can identify the sufficient statistic $T$ by simply factorizing the log-likelihood function. Hence no prior guess for T is needed. If $\theta_{\text{MLE}}$ denotes the maximum likelihood estimator (MLE) of $\theta$, then from \eqref{2eqn:factorization_thm}, we have
\begin{eqnarray*}
\theta_{\text{MLE}}
&=& \arg\max_\theta L({X}_1^n;\theta)\\
&=& \displaystyle\arg\max_\theta [u (\theta, T({X}_1^n)) + 
 v({X}_1^n)]\\
&=& \displaystyle\arg\max_\theta u (\theta, T({X}_1^n)).
\end{eqnarray*}
This implies that the MLE of $\theta$ depends on the sample only through the statistic $T$ (that is, knowledge of $T$ is sufficient for the estimation of $\theta$). This seems to suggest that the classical notion of sufficiency is based on ML estimation.

Observe that the `entire sample' is trivially a sufficient statistic with respect to any estimation. However, the notion of sufficient statistics is helpful only when it exists in a fixed number, independent of the sample size. \emph{Fisher-Darmois-Koopman-Pitman theorem} addresses this problem \cite{Fisher34J, Koopman36J, Pitman36J, Darmois35J, BarankinM63J}. It states that, under certain regularity conditions, the existence of such a set of sufficient statistics is possible if and only if the underlying family is exponential. A set of probability distributions $\mathcal{E}$ is said to be an \emph{exponential family} if there exist functions $w = [w_1,\dots,w_s]^\top$ and $f =$ $[f_1,\dots,$ $f_s]^\top$ with $w_i: \Theta \to \mathbb{R}$, $f_i: \mathbb{S}\to \mathbb{R}$, and $h:\mathbb{R}\to \mathbb{R}$ such that  each $p_\theta\in\mathcal{E}$ is given by  
	\begin{eqnarray}
	\label{2eqn:form_of_exp_family}
	{p_\theta(x)} = \left\{
	\begin{array}{ll}
	{Z(\theta)\exp\big[ h(x) + 
		w(\theta)^\top f(x) \big]} &\hbox{~if~} 
	{x}\in \mathbb{S}\\
	{0} &\hbox{~otherwise},
	\end{array}
	\right.
	\end{eqnarray}
where $Z(\theta)$ is the normalizing factor \cite{Brown86B}. Bernoulli, Binomial, Poisson, and Normal distributions are some examples of exponential families. However, a statistic that is sufficient for ML estimation might not be sufficient for some other estimations. For example, consider the family of binomial distributions with parameters, say $(m,\theta)$, where $m$ is a positive integer denoting the number of trials (assumed known) and $\theta\in(0,1)$ is the success probability, which is unknown. It is well-known that the sum of the sample is always a sufficient statistic for $\theta$ (with respect to the usual likelihood function) (see, for example, \cite{CasellaB02B}).

Let us now suppose that we want to estimate $\theta$ by maximizing the following likelihood function, called \emph{Cauchy-Schwarz likelihood function} (c.f. \cite[Eq. (2.90)]{Principe10B}):
\begin{eqnarray}
\label{2eqn:cauchyschwarz_distance}
L_{cs}(x_1^n;\theta)
\coloneqq  \log \Bigg[~{\sum\limits_{x=0}^m p_n(x) p_\theta(x)}\Big/{\sum\limits_{x=0}^m p_\theta(x)^2}\Bigg],
\end{eqnarray}
where $p_\theta(x) = \binom{m}{x}\theta^x (1-\theta)^{m-x}$ for $x=0,\dots,m$ and $p_n$ is the empirical distribution of the observed sample. 
Maximizing the Cauchy-Schwarz likelihood function is the same as minimizing the well-known Cauchy-Schwarz divergence. (This divergence is popular in machine learning applications such as blind source separation, independent component analysis, manifold learning, image segmentation, and so on \cite{JenssenPEE06J, ShanqingKW13J, KampaHP11conf, HoangVVM14J, NielsenSM17conf, NielsenSM17J}.)  
It is easy to check that the solutions of the estimating equation, in this case, are the roots of a polynomial in $\theta$ whose coefficients involve some of $p_n(x)$, $x=0,\dots,m$.
For example, when $m=2$, it is a polynomial of degree five, and the following are, respectively, the coefficients in descending order of exponents of $\theta$:\\

\noindent
$8p_n(1),~ 8p_n(0) -20p_n(1)-6p_n(2), ~ -20p_n(0)+16p_n(1)+12p_n(2), 18p_n(0)-4p_n(1)-6p_n(2),$
$-7p_n(0)-2p_n(1)+p_n(2)$, and $p_n(0)+p_n(1)$.
\vspace{0.2cm}

\noindent
This shows that at least two among $p_n(0), p_n(1)$ and $p_n(2)$ are needed to find the estimator. Hence the sum of the sample is no longer sufficient for estimating $\theta$. 

Observe that the Cauchy-Schwarz likelihood function is a special case ($\alpha=2$) of the more general class of likelihood functions, called {\em Jones et al. likelihood function} \cite{JonesHHB01J}:
 \begin{eqnarray}
\label{2eqn:likelihood_function_for_I_alpha}
{L_J^{(\alpha)}(x_1^n;\theta)}
&\coloneqq& \dfrac{1}{\alpha-1}\log\left[\dfrac{1}{n}\sum\limits_{i=1}^n p_\theta({x}_i)^{\alpha-1}\right] - \dfrac{1}{\alpha}\log\left[\int p_\theta(x)^\alpha
dx\right]\nonumber\\
&=& \dfrac{1}{\alpha-1}\log\left[\dfrac{1}{n}\sum\limits_{i=1}^n p_\theta({x}_i)^{\alpha-1}\right] - \log \|p_\theta\|,~~
 \end{eqnarray}
where $\alpha>0$, $\alpha\neq 1$, and $\|p_\theta\|\coloneqq \big(\int p_\theta(x)^\alpha dx\big)^{\frac{1}{\alpha}}$.

Let us now consider the following family of Student distributions given by
\begin{equation}
\label{2eqn:student_distribution_alpha2}
		p_{\theta}(x) = \frac{2}{\sqrt{3}\pi\sigma} \Big[1+ \frac{(x-\mu)^2}{3\sigma^2}\Big]^{-2}\quad\text{for}~x\in\mathbb{R},
\end{equation}
where $\theta=[\mu,\sigma^2]^\top$, $\mu\in (-\infty,\infty)$ and $\sigma\in(0,\infty)$.  Observe that (\ref{2eqn:student_distribution_alpha2}) does not fall in the exponential family as defined in (\ref{2eqn:form_of_exp_family}). Hence, by the Fisher-Darmois-Koopman-Pitman theorem, there does not exist any non-trivial sufficient statistics for $\mu$ and $\sigma$. Indeed, it can be shown that the knowledge of the entire sample is needed for estimating $\mu$ and $\sigma$ in this case (See \cite{GayenK21ISIT}). However, Eguchi et al. considered the estimation based on Jones et al. likelihood function (\ref{2eqn:likelihood_function_for_I_alpha}) with $\alpha =1/2$ and found closed-form solution for the estimators (See \cite{EguchiKK11J}). The estimators are functions of the statistics $\sum_{i=1}^nX_i$ and $\sum_{i=1}^nX_i^2$ \cite[Th. 3]{EguchiKK11J} (c.f. \cite[Th. 11]{GayenK21J}).
This motivates us to generalize the concept of sufficiency based on an estimation method.
Let us consider the difference of Jones et al. likelihood function $L_{J}^{(1/2)}$ evaluated at samples $x_1^n$ and $y_1^n$:
\begin{eqnarray*}
{L_{J}^{(1/2)}(x_1^n;\theta) - L_{J}^{(1/2)}(y_1^n;\theta)}
&=& 2\Big[\log \sum\limits_{i=1}^n p_\theta(y_i) - \log\sum\limits_{i=1}^n p_\theta(x_i)\Big].
\end{eqnarray*} 
For the family of Student distributions (\ref{2eqn:student_distribution_alpha2}), we have
\begin{eqnarray}
\label{2eqn:likelihood_difference_student_distribution}
{L_{J}^{(1/2)}(x_1^n;\theta) - L_{J}^{(1/2)}(y_1^n;\theta)}
&=& 2\log \Bigg[\dfrac{1+(1/(\mu^2 + 3\sigma^2))}{1+(1/(\mu^2+3\sigma^2))} \cdot  \dfrac{\sum\limits_{i=1}^n y_i^2 - (2\mu/(\mu^2+3\sigma^2)) \sum\limits_{i=1}^n y_i}{ \sum\limits_{i=1}^n x_i^2 - (2\mu/(\mu^2+3\sigma^2)) \sum\limits_{i=1}^n x_i}\Bigg].~
\end{eqnarray} 
(\ref{2eqn:likelihood_difference_student_distribution}) is independent of $\theta$ if
\begin{equation*}
    \sum\limits_{i=1}^n x_i^2 = \sum\limits_{i=1}^n y_i^2\quad\text{and}\quad
\sum\limits_{i=1}^n x_i = \sum\limits_{i=1}^n y_i.
\end{equation*}
Notice that the Jones et al. estimators for $\mu$ and $\sigma$ depend on the sample only through the statistics $\sum_{i=1}^n X_i$ and $\sum_{i=1}^n X_i^2$.
This implies that, when the estimation is based on Jones et al. likelihood function $L_{J}^{(1/2)}$, it is reasonable to call the pair $\big(\sum X_i, \sum X_i^2\big)$ sufficient for $\mu$ and $\sigma$.

 Jones et al. likelihood function is popular in the context of robust estimation \cite{JonesHHB01J}. Observe that (\ref{2eqn:likelihood_function_for_I_alpha}) coincides with the usual log-likelihood function $L(\theta)$ as $\alpha \to 1$. Hence Jones et al. likelihood function is considered a robust alternative to ML estimation. Motivated by the works of Windham \cite{Windham95J} and Field and Smith \cite{FieldS94J}, Jones et al. proposed this method in search of an estimating equation, an alternative to the one of MLE, that uses a suitable power of the model density to down-weight the effect of \emph{outliers} (possibly) present in the data \cite{JonesHHB01J, GayenK21J, Fujisawa13J, EguchiKK11J}. It is well-known that ML estimation is closely related to the minimization of the well-known Kullback-Leibler divergence or relative entropy \cite[Lem. 3.1]{CsiszarS04B}. Similarly, Jones et al. estimation is related to the minimization of the so-called Jones et al. divergence \cite{FujisawaE08J, LutwakYZ05J, Sundaresan02ISIT, Sundaresan07J, GayenK21J} (also known as $\gamma$-divergence \cite{CichockiA10J, FujisawaE08J}, projective power divergence \cite{EguchiKK11J}, logarithmic density power divergence \cite{MajiGB16J}, relative $\alpha$-entropy \cite{KumarS15J1, KumarS15J2}, R\'enyi pseudo-distance \cite{BroniatowskiTV12J, BroniatowskiV12J}). While Jones et al. divergence for $\alpha>1$ corresponds to robust inference, the same for $\alpha<1$ has connections to problems in information theory, such as mismatched guessing \cite{Sundaresan02ISIT}, source coding \cite{KumarS15J2}, encoding of tasks \cite{BunteL14J} and so on. Though Jones et al. estimation for $\alpha<1$ does not have a meaning in terms of robustness, it can be used to approximate MLE of a Student distribution (c.f. \cite[Rem. 7]{GayenK21J}).

Another likelihood function that is popular in the context of robust inference is the so-called Basu et al. likelihood function given by
\begin{eqnarray}
    \label{basu_likelihood}
    {L^{(\alpha)}_B(x_1^n;\theta)}
    &=& \frac{1}{n} \sum\limits_{i=1}^n
    \left[\frac{\alpha\cdot p_\theta(x_i)^{\alpha-1} -1}{\alpha -1}\right] -\int p_\theta(x)^{\alpha} dx,~~~
\end{eqnarray}
where $\alpha>0$, $\alpha\neq 1$ \cite[Eq. (2.2)]{BasuHHJ98J}. Maximization of $L^{(\alpha)}_B$ results in a non-normalized version of the Jones et al. estimating equation (see \cite[Sec.~1]{GayenK21J}). Basu et al. estimation is closely related to the minimization of the so-called Basu et al. divergence, which belongs to the Bregman class \cite{Csiszar95J1}. Basu et al. divergence for $\alpha = 2$ corresponds to the squared Euclidean distance between the empirical measure and the model density. A comparison of Basu et al. and Jones et al. estimations was studied in \cite{JonesHHB01J}. Unlike other likelihood functions such as the Hellinger likelihood function, Basu et al. and Jones et al. likelihood functions are special in robust inference, in the sense that smoothing of the empirical probability measure is not needed (see, for example, \cite[Sec.~1]{GayenK21J} for details).
In this paper, we propose a generalized notion of sufficiency based on the likelihood function associated with an estimation problem. 
Specifically, we will focus on Jones et al. and Basu et al. likelihood functions. We will characterize the family of probability distributions that has a fixed number of sufficient statistics independent of the sample size with respect to these likelihood functions. It turns out that these are the power-law families studied in \cite{GayenK21J, KumarS15J2, Bashkirov04J, Naudts04J, CsiszarM12J}. Student distributions fall in this family.

Sufficient statistics not only help in data reduction but also help in finding an estimator better than the given one. Indeed, if one has an unbiased estimator along with sufficient statistics, it is always possible to find another unbiased estimator (usually a function of sufficient statistics) whose variance is not more than the given one. This result is known as the \emph{Rao-Blackwell theorem} \cite{Rao45J, Blackwell47J}. That is, conditional expectation of an unbiased estimator given the sufficient statistics results in a uniform improvement (in variance) of the estimator. This suggests one to consider only such conditional expectations while searching for best unbiased estimators.

In this paper, we find Rao-Blackwell type estimators for power-law families of the form (\ref{2eqn:student_distribution_alpha2}) with respect to Jones et al. and Basu et al. likelihood functions. We also find a lower bound for the variance of these estimators. This generalizes the usual Cram\'er-Rao Lower Bound (CRLB) and results in a sharper bound for the power-law families. We also apply these results to the Student distributions.
\vspace{0.4cm}

Our main contributions in this paper are the following.
\vspace{0.2cm}

\begin{itemize}
	\item [(i)] Propose a generalized notion of principle of sufficiency that is appropriate when the underlying estimation is not necessarily ML estimation.
	 \vspace{0.2cm}
	
	\item [(ii)] Characterize the families of probability distributions that have a fixed number of sufficient statistics independent of sample size with respect to Basu et al. and Jones et al. likelihood functions. 
	 \vspace{0.2cm}
	 
	\item[(iii)] Generalize the notion of minimal sufficiency with respect to a likelihood function and find generalized minimal sufficient statistics for the power-law family.
	 \vspace{0.2cm}
	
	\item[(iv)] Extend Rao-Blackwell theorem to both Basu et al. and Jones et al. estimations and find best estimators for the power-law family. 
	 \vspace{0.2cm}
	
	\item[(v)] Derive a Cram\'er-Rao lower bound for the variance of an unbiased estimator (of the expected value of the generalized sufficient statistic) of a power-law family in terms of the asymptotic variance of the Jones et al. estimator.
 
\vspace{0.4cm}

\end{itemize}
The rest of the paper is organized as follows. In Section \ref{2sec:gen_notion_prin_suffi}, we introduce the generalized notion of sufficiency. In Section \ref{2sec:suff_stat_jones}, we establish Fisher-Darmois-Koopman-Pitman theorem for Jones et al. and Basu et al. likelihood functions. In this section, we also find minimal sufficient statistics for the mean and variance parameters of Student distributions. In Section \ref{2sec:Rao_blackwelLones_et_al}, we find generalized Rao-Blackwell type estimators for the power-law family. In Section \ref{2sec:application_student_dist}, we establish generalized CRLB for Basu et al. and Jones et al. estimators of power-law family. In this section, we also find the best estimators for the mean and variance parameters of Student distributions. We end the paper with a summary and concluding remarks in Section \ref{2sec:summary}. Proofs of most of the results are presented in the Appendix.

\section{A Generalized notion of the principle of sufficiency}
\label{2sec:gen_notion_prin_suffi}
In this section, we first propose a generalized notion of sufficiency analogous to Koopman's definition of sufficiency. We then derive a factorization theorem associated with this generalized notion. We also show that this is equivalent to generalized Fisher's definition of sufficiency. 
\vspace{0.2cm}

\begin{definition}
\label{1defn:gen_koopman}
  Let $\Pi=\{p_\theta:\theta\in\Theta\}$ be a family of probability distributions. Suppose the underlying problem is to estimate $\theta$ by maximizing some generalized likelihood function $L_G$.
A set of statistics $T$ is said to be a {\em sufficient statistic for $\theta$ with respect to $L_G$} if $[L_G({x}_1^n;\theta)-L_G({y}_1^n;\theta)]$ is independent of $\theta$ whenever $T({x}_1^n) = T({y}_1^n)$. 
\end{definition}
\vspace{0.2cm}

Observe that when $L_G$ is the log-likelihood function, the above definition coincides with the usual Koopman's definition of sufficiency. We now derive a generalized factorization theorem.
\vspace{0.2cm}

\begin{proposition}[Generalized Factorization Theorem]
	\label{2prop:gen_factorization_thm}
	Let $\Pi=\{p_\theta:\theta\in\Theta\}$ be a family of probability distributions.  Suppose the underlying problem is to estimate $\theta$ by maximizing some generalized likelihood function $L_G$. Then a set of statistics $T$ is a sufficient statistic for $\theta$ with respect to $L_G$ if and only if there exist functions $u:\Theta\times\mathscr{T}\to\mathbb{R}$ and $v:\mathbb{S}^n\to\mathbb{R}$
	such that
	\begin{equation}
\label{2eqn:Gen_Factor_Thm}
	L_G({x}_1^n;\theta) = u (\theta, T({x}_1^n)) + v({x}_1^n),
\end{equation} 
	for all sample points ${x}_1^n$ and for all $\theta\in \Theta$.
\end{proposition}
\vspace{0.2cm}

\begin{IEEEproof}
   See Sec-A.1 of Appendix.
\end{IEEEproof}
\vspace{0.2cm}

Observe that when $L_G=L$, the above coincides with the usual factorization theorem \cite{Neyman35J, HalmosS49J}.
We now show that this generalized notion of sufficiency is also equivalent to a generalized version of Fisher's definition of sufficiency. Let $X_1^n$ be i.i.d. random variables according to some $p_\theta\in\Pi$. Let $x_1^n\coloneqq x_1,\dots,x_n$ be a sample. Let us define the deformed probability distribution associated with the likelihood function $L_G$ by
\begin{eqnarray*}
    p_\theta^*(x_1^n) \coloneqq \frac{\exp[ L_G(x_1^n ; \theta)]}{\int \exp[ L_G(y_1^n ; \theta)] dy_1^n},
\end{eqnarray*}
provided the denominator is finite. Observe that, if $L_G(x_1^n;\theta) = L(x_1^n;\theta)$, then
\begin{eqnarray}
    \label{2eq:p_star_coincies_p_jones_et_al}
   p_\theta^*(x_1^n) &=& \dfrac{\exp[\log p_\theta(x_1^n)]}{\int \exp[\log p_\theta(y_1^n)] d{y_1^n}}
   = \dfrac{ p_\theta(x_1^n)}{\int p_\theta(y_1^n) d{y_1^n}}
   = p_\theta(x_1^n).
\end{eqnarray}

Let $T$ be a statistic of the sample $x_1^n$. For any value $t$ of $T$, let $q_\theta^*(t)\coloneqq\int p_\theta^*(y_1^n)\textbf{1}(T(y_1^n)=t) dy_1^n$, where $\textbf{1}(\cdot)$ denotes the indicator function. This can be re-written as $$q_\theta^*(t)= \int\limits_{A_t} p_\theta^*(y_1^n) dy_1^n,$$ where $A_t=\{x_1^n: T(x_1^n)=t\}$.
The associated conditional distribution is given by
\begin{eqnarray}
	\label{2eqn:gen_condition_dist}
	p_{\theta{_{X_1^n|T}}}^*(x_1^n|t) &\coloneqq & \dfrac{p_\theta^*(x_1^n)\textbf{1}(T(x_1^n) = t)}
      {q_\theta^*(t)} \nonumber\\
      &=& 
      \left\{
	\begin{array}{ll}
	{\dfrac{p_\theta^*(x_1^n)}{q_\theta^*(t)}} &\hbox{~if~} 
	x_1^n\in A_t\\
	{0} &\hbox{~otherwise}
	.
	\end{array}
	\right.
\end{eqnarray}
Then
\begin{eqnarray*}
\int p_{\theta{_{X_1^n|T}}}^*(x_1^n|t) d{x_1^n}
= \int\limits_{A_t} \dfrac{p_\theta^*(x_1^n)}{q_\theta^*(t)} dx_1^n
= \int\limits_{A_t} \dfrac{p_\theta^*(x_1^n)}{q_\theta^*(T(x_1^n))} dx_1^n
= \dfrac{\int\limits_{A_t} p_\theta^*(x_1^n) dx_1^n}{\int\limits_{A_t} p_\theta^*(x_1^n)dx_1^n}
= 1.
\end{eqnarray*}
Observe that, when $L_G = L$, we have
\begin{equation*}
 p_{\theta{_{X_1^n|T}}}^*(x_1^n|t) =  p_{\theta{_{X_1^n|T}}}(x_1^n|t),
\end{equation*}
the usual conditional distribution of $x_1^n$ given $T=t$.
\vspace{0.2cm}

\begin{proposition}
   \label{equi_generalized_fisher_generalized_facto}
$T$ is sufficient for $\theta$ with respect to a likelihood function $L_G$ if and only if $p_{\theta{_{X_1^n|T}}}^*(x_1^n|t)$ is independent of $\theta$ for any value $t$ of $T$. 
\end{proposition}
\vspace{0.2cm}

\begin{IEEEproof}
   See Sec-A.2 of Appendix.
\end{IEEEproof}
\vspace{0.2cm}

\begin{remark}
Fisher’s idea of sufficiency has one appealing interpretation that if we know the value of the sufficient statistic $T$ without knowing the entire sample, then we could (using the conditional distribution of the sample which is independent of $\theta$ given the sufficient statistic) generate another sample which has the same unconditional distribution as the original sample and the same value of the sufficient statistic $T$. Hence the information contained in both the samples are essentially the same. In this sense, the statistic $T$ is ``sufficient''. This interpretation continues to hold for the generalized notion of sufficiency as well. Indeed, in the case of generalized sufficient statistics $T$ with respect to $L_G$, one can generate a sample using the generalized conditional distribution as defined in \eqref{2eqn:gen_condition_dist}. It contains the same information about $\theta$ as in the original sample  (with respect to $L_G$). Suppose that the given sample is $X_1^n=x_1^n$ and that $T(x_1^n) = t$. If we do not know $x_1^n$, but know only that $T(X_1^n)=t$, then we can generate $z_1^n$ using the generalized conditional distribution of $X_1^n=z_1^n$, given $T(X_1^n)=t$ (since the conditional distribution is independent of $\theta$ by Proposition \ref{equi_generalized_fisher_generalized_facto}). Then, using \eqref{2eqn:like_fact_expre}, the value of $L_G$ for the generated sample is given by $L_G(z_1^n;\theta) = u(t,\theta) +v(z_1^n)$, where 
\begin{eqnarray*}
    u(t,\theta) &=& \log \int\limits_{A_t} \exp[L_G(y_1^n;\theta)] dy_1^n\quad\text{with}\\
    A_t&=&\{x_1^n: T(x_1^n)=t\}, \text{ and }\\\\
 v(z_1^n)& =& \log p_{\theta_{X_1^n|T}}^*(z_1^n|t).
\end{eqnarray*}
Hence
\begin{eqnarray*}
    \arg\displaystyle\max_\theta L_G(z_1^n;\theta) = \arg\displaystyle\max_\theta u(t,\theta)
    =  \arg\displaystyle\max_\theta L_G(x_1^n;\theta),
\end{eqnarray*}
since $x_1^n\in A_t$. Thus $x_1^n$ and $z_1^n$ have the same information about $\theta$ as far as estimation with respect to $L_G$ is concerned. Also, from \eqref{2eqn:gen_condition_dist}, for any $y_1^n\in A_t$, the unconditional distribution with respect to $p_\theta^*$ is given by $p^*_{\theta}(y_1^n) = p^*_{\theta_{X_1^n|T}}(y_1^n|t) q_\theta^*(t)$. This implies that the generated sample $z_1^n$ has the same distribution as the initial sample $x_1^n$ (with respect to $p_\theta^*$). 
\end{remark}
\vspace{0.2cm}

While sufficient statistics help in data compression, \emph{minimal sufficient statistics} help maximum such compression. We now define the notion of minimal sufficiency in a general setup. Let us first recall the definition of minimal sufficiency in the classical sense.
\vspace{0.2cm}

\begin{definition}\cite[Def. 6.2.11]{CasellaB02B}
	\label{1defn:minimal_sufficient_statistics}
	A statistic $T$ is said to be a \emph {minimal sufficient statistic} if $T$ is a function of any other sufficient statistic. In other words, $T$
	is minimal if the following holds: For any sufficient statistic $\widetilde{T}$ and for any two i.i.d. samples
	${x}_1^n$ and ${y}_1^n$, $\widetilde{T}(x_1^n) =
	\widetilde{T}(y_1^n)$ implies 
	${T}({{x}_1^n}) =
	{T}({y}_1^n)$.
\end{definition}
\vspace{0.2cm}

There is an easy criterion to find minimal sufficient statistics for MLE due to Lehmann and Scheff\'e \cite[Th. 6.2.13]{CasellaB02B}. Here we first generalize this criterion to a generalized likelihood function. Later we will use it to find minimal sufficient statistics for power-law families.
\vspace{0.2cm}

\begin{proposition}
	\label{1prop:condition_for_minimal_sufficient_statistic}
	Let $\Pi=\{p_\theta:\theta\in\Theta\}$ be a family of probability distributions. Suppose the underlying problem is to estimate $\theta$ by maximizing some likelihood function $L_G$. Assume that the following condition holds for a statistic $T$. For any two i.i.d. samples
	$x_1^n$ and $y_1^n$,
	$T(x_1^n) = T(y_1^n)$ if and only if
	$[L_G(x_1^n;\theta) - L_G(y_1^n;\theta)]$ is independent of $\theta$ for $\theta \in \Theta$. Then $T$ is a minimal sufficient statistic for $\theta$ with respect to $L_G$.
\end{proposition}
\vspace{0.2cm}

\begin{IEEEproof}
	See Sec-A.3 of Appendix.
\end{IEEEproof}

\section{Generalized Fisher-Darmois-Koopman-Pitman Theorem}
\label{2sec:suff_stat_jones}
In this section, we first characterize the probability distributions that have a fixed number of sufficient statistics independent of sample size with respect to Jones et al. and Basu et al. likelihood functions. We observe that each of them is a power-law extension of the usual exponential family. We then find sufficient statistics for the power-law families. We start by defining the family associated with Jones et al. likelihood function.

\begin{definition}
	\label{2defn:general_M_alpha_family}
	A {\em $k$-parameter $\alpha$-power-law family} characterized by $h,w,f,\Theta$ and $\mathbb{S}$, where $\Theta$ is an open subset of $\mathbb{R}^k$, $w = [w_1,\dots,w_s]^\top$, $f =$ $[f_1,\dots,$ $f_s]^\top$, with $w_i: \Theta \to \mathbb{R}$ is differentiable for $i=1,\dots,s$, $f_i: \mathbb{R}\to \mathbb{R}$ for $i=1,\dots,s$, $h:\mathbb{R}\to \mathbb{R}$, and $\mathbb{S}\subseteq \mathbb{R}$, is defined by $\mathbb{M}^{(\alpha)}=\{p_\theta\colon \theta\in\Theta\}$, where
	\begin{eqnarray}
	\label{2eqn:form_of_general_M_alpha_family}
	{p_\theta(x)} = \left\{
	\begin{array}{ll}
	{Z(\theta)\big[ h(x) +  
		w(\theta)^\top f(x) \big]^{\frac{1}{\alpha - 1}}} &\hbox{~for~} 
	x \in \mathbb{S}\\
	{0} &\hbox{~otherwise,}
	\end{array}
	\right.
	\end{eqnarray}
	with $Z(\theta)^{-1} = \int\limits_{\mathbb{S}}[h(x) +  w(\theta)^\top f(x)]^\frac{1}{\alpha-1} d{x}$.
\end{definition}
\vspace{0.2cm}

\begin{remark}
\label{2rem:remark_on_m_alpha_family}
\begin{itemize}
    \item [(i)] If $h(x)=q(x)^{\alpha-1}$ for some function $q(x)>0$ and $f(x)=(1-\alpha)\Tilde{f}(x)$, then $\mathbb{M}^{(\alpha)}$-family can be shown to coincide with an exponential family when $\alpha = 1$ (see \cite[Def.~8]{KumarS15J1}).

\vspace{0.2cm}

\item[(ii)] In \cite{KumarS15J2}, this family was derived as projections of Jones et al. divergence on a set of probability distributions determined by linear constraints. More about $\mathbb{M}^{(\alpha)}$-family can be found in \cite{Bashkirov04J, GayenK18ISIT, GayenK21J}, and \cite{KumarS15J2}.
\vspace{0.2cm}

\item[(iii)] Student distributions fall in $\mathbb{M}^{(\alpha)}$-family \cite[Ex.~3]{GayenK21J}.
\end{itemize}
\end{remark}
\vspace{0.2cm}

We now state the main result of this section.
\vspace{0.2cm}

\begin{theorem}
	\label{2thm:necessary_condition_sufficiency_jones}
	Let $\Pi=\{p_\theta:\theta\in\Theta\subset \mathbb{R}^k\}$ be a family of probability distributions with common support $\mathbb{S}$, an open subset of $\mathbb{R}$.
	Let $x_1^n$ be an i.i.d. sample from some $p_\theta\in\Pi$. For some $s$ $(k\leq s<n)$, let  $T_j(x_1^n)$ for $j=1,\dots,s$ be differentiable functions on $\mathbb{S}^n$.
	 If $[T_1,\dots,T_s]^\top$ is a sufficient statistic for $\theta$ with respect to Jones et al. likelihood function $L_J^{(\alpha)}$, then there exist functions $h:\mathbb{S}\to  \mathbb{R}$, $f_i:\mathbb{S}\to \mathbb{R}$ and $w_i:\Theta\to \mathbb{R}$,  $i=1,\dots,s$ such that $p_\theta(x)$ can be expressed as (\ref{2eqn:form_of_general_M_alpha_family}) for $x\in \mathbb{S}$.  
\end{theorem}
\vspace{0.2cm}

\begin{IEEEproof}
	The proof is divided into two parts. The first part deals with the case when $k=s$ and the next one for $k< s$.
\begin{itemize}
    \item [(a)] Let $k=s$.
	We only give the proof when $k=s=1$. Other cases can be proved similarly. 
	
	Let us first consider the Jones et al. likelihood function (\ref{2eqn:likelihood_function_for_I_alpha}).
Taking derivative with respect to $\theta$ on both sides of (\ref{2eqn:likelihood_function_for_I_alpha}), we get
\begin{eqnarray}
    \label{2eqn:likelihood_derivative_jones}
{\frac{\partial}{\partial\theta} L_J^{(\alpha)}(x_1^n;\theta)} 
&=&
\frac{1}{\alpha -1} \frac{\sum_{i=1}^n \frac{\partial}{\partial\theta} \big[p_\theta (x_i)^{\alpha -1}\big]}{\sum_{i=1}^n p_\theta(x_i)^{\alpha -1}} - \frac{\partial}{\partial\theta} \log \|p_\theta\|.\nonumber\\
\end{eqnarray}
If $T$ is a sufficient statistic for $\theta$ with respect to $L_J^{(\alpha)}$, then, from Proposition \ref{2prop:gen_factorization_thm}, we have
\begin{equation}
\label{2eqn:factorization_of_jones_likelihood}
L_J^{(\alpha)}(x_1^n;\theta) = u(\theta, T) + v(x_1^n)
\end{equation}
for some functions $u$ and $v$, where $v$ is independent of $\theta$. This implies that
\begin{equation}
\label{2eqn:derivative_facorization_jones}
\dfrac{\partial }{\partial\theta}L_J^{(\alpha)}(x_1^n;\theta) = u_1(\theta, T),
\end{equation}
where $u_1(\theta, T) \coloneqq \dfrac{\partial}{\partial\theta}\big[u(\theta, T)\big]$. Comparing (\ref{2eqn:likelihood_derivative_jones}) and (\ref{2eqn:derivative_facorization_jones}), we have
\begin{eqnarray}
    \label{2eqn:likelihood_and_factorization_derivative_jones}
{u_1(\theta, T)} 
&=& \frac{1}{\alpha -1} \frac{\sum_{i=1}^n \frac{\partial}{\partial\theta} \big[p_\theta (x_i)^{\alpha -1}\big]}{\sum_{i=1}^n p_\theta(x_i)^{\alpha -1}} -\frac{\partial}{\partial\theta} \log \|p_\theta\|.\nonumber\\
\end{eqnarray}
Now, for a particular value of $\theta\in\Theta$, say $\theta_0$, (\ref{2eqn:likelihood_and_factorization_derivative_jones}) can be re-written as
\begin{align}
\label{2eqn:likelihood_and_factorization_derivative_jones_1}
\lefteqn{u_1(\theta_0, T)}\nonumber\\
& = \frac{1}{\alpha -1} \frac{\sum_{i=1}^n \Big[\frac{\partial}{\partial\theta} \big\{p_\theta (x_i)^{\alpha -1}\big\}\Big]_{\theta_0}}{\sum_{i=1}^n p_{\theta_0}(x_i)^{\alpha -1}}- \Big[\frac{\partial}{\partial\theta} \log \|p_\theta\|\Big]_{\theta_0}\nonumber\\
&= \frac{1}{\alpha -1} \frac{\sum_{i=1}^n g(x_i)}{\sum_{i=1}^n h(x_i)} - \Big[\frac{\partial}{\partial\theta} \log\|p_\theta\|\Big]_{\theta_0},
\end{align}
where $g(x_i) = \big[\frac{\partial}{\partial\theta} \big\{p_\theta (x_i)^{\alpha -1}\big\}\big]_{\theta_0}$, and $h(x_i) = p_{\theta_0}(x_i)^{\alpha -1}$.

Let us denote $\widetilde{T}_1\coloneqq \sum_{i=1}^n g(x_i)$ and $\widetilde{T}_2\coloneqq \sum_{i=1}^n h(x_i)$.
Then, from (\ref{2eqn:likelihood_and_factorization_derivative_jones_1}), there exists a function $\Phi$ such that
\begin{equation}
\label{2eqn:form_of_T}
T = \Phi(\widetilde{T}_1,\widetilde{T}_2).
\end{equation}
Using this,
(\ref{2eqn:derivative_facorization_jones}) can be re-written as
\begin{equation}
\label{2eqn:likelihood_two_suff_stat_jones}
\dfrac{\partial }{\partial\theta}L_J^{(\alpha)}(x_1^n;\theta) = u_2(\theta, \widetilde{T}_1,\widetilde{T}_2)
\end{equation}
for some function $u_2$. Observe also that, (\ref{2eqn:likelihood_function_for_I_alpha}) can be re-written as
\begin{align*}
\lefteqn{L_J^{(\alpha)}(x_1^n;\theta)}\\
&= \dfrac{1}{\alpha-1}\log\bigg[\dfrac{1}{n}\sum\limits_{i=1}^n p_\theta({x}_i)^{\alpha-1}\bigg]
- \dfrac{1}{\alpha-1} \log \|p_\theta\|^{\alpha -1}\\
& = \dfrac{1}{\alpha-1}\log \Bigg[\frac{1}{n}\sum_{i=1}^n \frac{p_\theta(x_i)^{\alpha -1}}{\|p_\theta\|^{\alpha -1}}\Bigg].
\end{align*}
This implies that
\begin{equation*}
\exp \Big[(\alpha -1) L_J^{(\alpha)}(x_1^n;\theta)\Big] =
\frac{1}{n}\sum_{i=1}^n \frac{p_\theta(x_i)^{\alpha -1}}{\|p_\theta\|^{\alpha -1}}.
\end{equation*}
Further (\ref{2eqn:factorization_of_jones_likelihood}) implies that
\begin{eqnarray*}
    {\exp \Big[(\alpha -1) L_J^{(\alpha)}(x_1^n;\theta)\Big]} 
    &=& 
\exp\Big[ (\alpha -1) u(T,\theta)\Big]\cdot \exp\Big[ (\alpha -1) v(x_1^n)\Big].
\end{eqnarray*}
Thus we have
\begin{eqnarray}
    \label{2eqn:exponential_of_likelihood_jones_factorization}
{\frac{1}{n}\sum_{i=1}^n \frac{p_\theta(x_i)^{\alpha -1}}{\|p_\theta\|^{\alpha -1}}}
&=&
\exp\big[ (\alpha -1) u(T,\theta)\big]\cdot \exp\big[ (\alpha -1) v(x_1^n)\big].\nonumber\\
\end{eqnarray}
Fixing $\theta=\theta_0$, the above equation becomes
\begin{eqnarray*}
   \widetilde{T}_2 &=& n\|p_{\theta_0}\|^{\alpha -1}\exp\Big[ (\alpha -1) u(T,\theta_0)\Big]\exp\Big[ (\alpha -1) v(x_1^n)\Big].
\end{eqnarray*}
Since $T=\Phi(\widetilde{T}_1,\widetilde{T}_2)$, this implies that $v(x_1^n)$ must be a function of $\widetilde{T}_1$ and $\widetilde{T}_2$ only. Therefore, (\ref{2eqn:exponential_of_likelihood_jones_factorization}) can be re-written as
\begin{equation}
\label{2eqn:distribution_from_likelihood_jones}
\sum_{i=1}^n {p_\theta(x_i)^{\alpha -1}} =
u_3(\theta,\widetilde{T}_1,\widetilde{T}_2),
\end{equation}
where $u_3(\theta, \widetilde{T}_1,\widetilde{T}_2)\coloneqq n\|p_\theta\|^{\alpha -1}\exp\big[(\alpha -1) u\big(\Phi(\widetilde{T}_1,\widetilde{T}_2),\theta\big)\big] \cdot\exp\big[(\alpha -1) v(x_1^n)\big]$.
Now, taking partial derivative with respect to $x_j$ on both sides of (\ref{2eqn:distribution_from_likelihood_jones}), we have
\begin{align}
\label{2eqn:derivative_p_theta_wrt_X_j}
\frac{\partial}{\partial x_j} p_\theta(x_j)^{\alpha -1} 
&= \frac{\partial}{\partial \widetilde{T}_1} u_3(\theta,\widetilde{T}_1,\widetilde{T}_2) \frac{\partial \widetilde{T}_1}{\partial x_j} + \frac{\partial}{\partial \widetilde{T}_2} u_3(\theta,\widetilde{T}_1,\widetilde{T}_2) \frac{\partial \widetilde{T}_2}{\partial x_j}\nonumber\\
& = \frac{\partial}{\partial \widetilde{T}_1} u_3(\theta,\widetilde{T}_1,\widetilde{T}_2) \frac{\partial }{\partial x_j}g(x_j) + \frac{\partial}{\partial \widetilde{T}_2} u_3(\theta,\widetilde{T}_1,\widetilde{T}_2) \frac{\partial}{\partial x_j}h(x_j).
\end{align}
Since the left hand side of (\ref{2eqn:derivative_p_theta_wrt_X_j}) is a function of $x_j$ and $\theta$ only, we shall denote it by $A(x_j,\theta)$. Thus taking partial derivative with respect to $x_i$ for $i\neq j$ on both sides of (\ref{2eqn:derivative_p_theta_wrt_X_j}), we get
\begin{eqnarray}
    \label{2eqn:derivative}
&&\hspace{-1cm}\frac{\partial}{\partial x_i}  \Big[\frac{\partial}{\partial \widetilde{T}_1} u_3(\theta, \widetilde{T}_1,\widetilde{T}_2)\Big]~  \frac{\partial }{\partial x_j}[g(x_j)] + \frac{\partial}{\partial x_i}  \Big[\frac{\partial}{\partial \widetilde{T}_2}u_3(\theta, \widetilde{T}_1,\widetilde{T}_2)\Big]~ \frac{\partial}{\partial x_j}[h(x_j)] = 0.~~~
\end{eqnarray}
From (\ref{2eqn:derivative}), it is clear that 
\begin{eqnarray}
    \label{2eqn:zero_derivative}
&&\hspace{-1cm}\frac{\partial}{\partial x_i}  \Big[\frac{\partial}{\partial \widetilde{T}_1} u_3(\theta, \widetilde{T}_1,\widetilde{T}_2)\Big] = 0~\text{if and only if}~\frac{\partial}{\partial x_i}  \Big[\frac{\partial}{\partial \widetilde{T}_2}u_3(\theta, \widetilde{T}_1,\widetilde{T}_2)\Big] =0,
\end{eqnarray}
since neither $\frac{\partial }{\partial x_j}[g(x_j)]$ nor $\frac{\partial }{\partial x_j}[h(x_j)]$ can be zero.
 
\noindent  
Equivalently, 
\begin{eqnarray*}
    &&\hspace{-1cm}\frac{\partial}{\partial x_i}  \Big[\frac{\partial}{\partial \widetilde{T}_1} u_3(\theta, \widetilde{T}_1,\widetilde{T}_2)\Big] \neq 0~\text{if and only if}~\frac{\partial}{\partial x_i}  \Big[\frac{\partial}{\partial \widetilde{T}_2}u_3(\theta, \widetilde{T}_1,\widetilde{T}_2)\Big] \neq 0.
\end{eqnarray*}

Suppose that $\dfrac{\partial}{\partial x_i}  \Big[\dfrac{\partial}{\partial \widetilde{T}_1} u_3(\theta, \widetilde{T}_1,\widetilde{T}_2)\Big] \neq 0$. Then (\ref{2eqn:derivative}) can be re-written as
\begin{eqnarray}
    \label{2eqn:ratio_derivatives}
\dfrac{\frac{\partial}{\partial x_i}  \Big[\frac{\partial}{\partial \widetilde{T}_1} u_3(\theta, \widetilde{T}_1,\widetilde{T}_2)\Big]}{\frac{\partial}{\partial x_i}  \Big[\frac{\partial}{\partial \widetilde{T}_2}u_3(\theta, \widetilde{T}_1,\widetilde{T}_2)\Big]} &=&
-\dfrac{\frac{\partial}{\partial x_j}[h(x_j)]}{\frac{\partial}{\partial x_j}[g(x_j)]}\eqqcolon \xi(x_j).
\end{eqnarray}
Since the right hand side of (\ref{2eqn:ratio_derivatives}) is a function of $x_j$ alone, both $\frac{\partial}{\partial \widetilde{T}_1} [u_3(\theta, \widetilde{T}_1,\widetilde{T}_2)]$ and $\frac{\partial}{\partial \widetilde{T}_2} [u_3(\theta, \widetilde{T}_1,\widetilde{T}_2)]$ must be of the following form
\begin{eqnarray*}
&&\hspace{-1.2cm} \frac{\partial}{\partial \widetilde{T}_1} u_3(\theta, \widetilde{T}_1,\widetilde{T}_2) = \Psi(x_1^n,\theta)S_1(x_j) + V_1(x_1,\dots,x_{i-1},x_{i+1},\dots,x_n,\theta)
\end{eqnarray*}
and
\begin{eqnarray*}
&&\hspace{-1.2cm}\frac{\partial}{\partial \widetilde{T}_2} u_3(\theta, \widetilde{T}_1,\widetilde{T}_2) = \Psi(x_1^n,\theta)S_2(x_j)+ V_2(x_1,\dots,x_{i-1},x_{i+1},\dots,x_n,\theta),
\end{eqnarray*}
for some real-valued functions $\Psi, S_1,S_2,V_1$, and $V_2$.
Observe that (\ref{2eqn:ratio_derivatives}) holds for any $i\neq j$. This implies $V_1$ and $V_2$ must be functions of $x_j$ and $\theta$ only.
Using these in (\ref{2eqn:derivative_p_theta_wrt_X_j}), we get
\begin{equation}
\label{2eqn:form_of_f}
A(x_j,\theta) = \Psi(x_1^n,\theta) K_1(x_j) + K_2(x_j,\theta),
\end{equation}
where 
\begin{eqnarray*}
{K_1(x_j)}
&\coloneqq& S_1(x_j)~ \frac{\partial }{\partial x_j}[g(x_j)] - S_2(x_j)~ \frac{\partial }{\partial x_j}[h(x_j)],
\end{eqnarray*}
and
\begin{eqnarray*}
{K_2(x_j,\theta)}
&\coloneqq& V_1(x_j,\theta)~ \frac{\partial }{\partial x_j}[g(x_j)] - V_2(x_j,\theta)~ \frac{\partial }{\partial x_j}[h(x_j)].
\end{eqnarray*}
\noindent
In order for (\ref{2eqn:form_of_f}) to hold, $\Psi$ must be a function of $x_j$ and $\theta$ only. This implies that both 
\begin{equation}
\label{2eqn:deriavtive_g}
\dfrac{\partial}{\partial \widetilde{T}_1} u_3(\theta, \widetilde{T}_1,\widetilde{T}_2)\quad\text{and}\quad \dfrac{\partial}{\partial \widetilde{T}_2} u_3(\theta, \widetilde{T}_1,\widetilde{T}_2)
\end{equation}
must be a function of $x_j$ and $\theta$ only. 

On the other hand, let $\frac{\partial}{\partial x_i}  \Big[\frac{\partial}{\partial \widetilde{T}_1} u_3(\theta, \widetilde{T}_1,\widetilde{T}_2)\Big] = 0$ for some $i\neq j$. Then from (\ref{2eqn:zero_derivative}), we have
$\frac{\partial}{\partial x_i}  \Big[\frac{\partial}{\partial \widetilde{T}_2} u_3(\theta, \widetilde{T}_1,\widetilde{T}_2)\Big] = 0$.
This implies that both $\frac{\partial}{\partial \widetilde{T}_1} u_3(\theta, \widetilde{T}_1,\widetilde{T}_2)$ and $\frac{\partial}{\partial \widetilde{T}_2} u_3(\theta, \widetilde{T}_1,\widetilde{T}_2)$ are independent of $x_i$'s where $i\neq j$. Thus, in this case also, (\ref{2eqn:deriavtive_g}) holds.
\vspace{0.1cm}

Since $u_3(\theta,\widetilde{T}_1,\widetilde{T}_2)$ is a function of $\widetilde{T}_1$, $\widetilde{T}_2$ and $\theta$ only, both $\frac{\partial}{\partial \widetilde{T}_1} u_3(\theta, \widetilde{T}_1,\widetilde{T}_2)$ and
$\frac{\partial}{\partial \widetilde{T}_2} u_3(\theta, \widetilde{T}_1,\widetilde{T}_2)$ are also functions of $\widetilde{T}_1,\widetilde{T}_2$ and $\theta$ only.
But both $\widetilde{T}_1$ and $\widetilde{T}_2$ are symmetric functions of the $x_i$'s, $i=1,\dots,n$.
This implies that both quantities in (\ref{2eqn:deriavtive_g}) must only be a function of $\theta$. Let 
\begin{eqnarray*}
  &&\hspace{-2cm}\frac{\partial}{\partial \widetilde{T}_1} u_3(\theta,\widetilde{T}_1,\widetilde{T}_2)=:w_1(\theta)\quad\text{and}
\quad\frac{\partial}{\partial \widetilde{T}_2} u_3(\theta,\widetilde{T}_1,\widetilde{T}_2)=:w_2(\theta)
\end{eqnarray*}
for some  $w_1$ and $w_2$ which are functions of $\theta$ only. Then we have
\begin{equation*}
u_3(\theta,\widetilde{T}_1,\widetilde{T}_2) = w_1(\theta) \widetilde{T}_1 + w_2(\theta) \widetilde{T}_2 +w_3(\theta).
\end{equation*}
Using (\ref{2eqn:distribution_from_likelihood_jones}), we get
\begin{equation*}
p_\theta(x_i)^{\alpha -1} = w_1(\theta) g(x_i) + w_2(\theta) h(x_i) +w_3(\theta)
\end{equation*} 
for $i=1,\dots,n$. 
\item[(b)] Let $k<s$. The proof for this case is quite similar to that of case (a). For the sake of simplicity, we give a sketch of proof for $k=1$ and $s=2$. 

Let $(T_1,T_2)$ be a pair of sufficient statistic of $\theta$ with respect to Jones et al. likelihood function $L_J^{(\alpha)}$. Then, proceeding as in Case (a), one would get an equation similar to \eqref{2eqn:form_of_T} with four symmetric functions of $x_i$'s (instead of two), say, $\widetilde{T}_1= \sum_{i=1}^n g_1(x_i)$, $\widetilde{T}_2= \sum_{i=1}^n g_2(x_i)$, $\widetilde{T}_3= \sum_{i=1}^n h_1(x_i)$ and $\widetilde{T}_4= \sum_{i=1}^n h_2(x_i)$. Hence \eqref{2eqn:distribution_from_likelihood_jones} holds with $\widetilde{T}_1,\widetilde{T}_2,\widetilde{T}_3$, and $\widetilde{T}_4$ instead of only $\widetilde{T}_1$ and $\widetilde{T}_2$.
Then a similar calculation yields that all $\frac{\partial}{\partial \widetilde{T}_i} u_3(\theta,\widetilde{T}_1,\widetilde{T}_2,\widetilde{T}_3,\widetilde{T}_4)$ for $i=1,2,3,4$ are a function of $\theta$ only. This further implies
\begin{eqnarray*}
    {p_\theta(x_i)^{\alpha -1} = \sum\limits_{j=1}^2 w_j^{(1)}(\theta) g_j(x_i)}  + \sum\limits_{j=1}^2 w_j^{(2)}(\theta) h_j(x_i) + w^{(3)}(\theta)
\end{eqnarray*}
for $i=1,\dots,n$. This completes the proof.
\end{itemize}
\end{IEEEproof}
\vspace{0.2cm}

The converse of the above theorem is also true. Indeed, in the following, we show that an $\mathbb{M}^{(\alpha)}$-family always has a fixed number of sufficient statistics with respect to the Jones et al. likelihood function $L_J^{(\alpha)}$. 
\vspace{0.2cm}

\begin{theorem}
	\label{2thm:suffiicent_statistics_M_alpha}
	Consider a parametric family of probability distributions $\Pi=\{p_\theta:\theta\in\Theta\subset \mathbb{R}^k\}$ having a common support $\mathbb{S}$.
	Let $T=[T_1,\dots,T_s]^\top$ be $s$ statistics defined on $\mathbb{S}^n$ where $s\geq k$.  $T$ is a sufficient statistic for $\theta$ with respect to Jones et al. likelihood function $L_J^{(\alpha)}$ if the following hold.
	\begin{itemize}
	    \item [(a)] $\Pi$ is a $k$-parameter $\mathbb{M}^{(\alpha)}$-family as in Definition \ref{2defn:general_M_alpha_family},
	    \item[(b)] there exist real valued
	functions $\psi_1,\dots,\psi_s$ defined on $\mathscr{T}$ such that for any i.i.d. sample $x_1^n$
	\begin{equation}
	\psi _i (T({x}_1^n)) = \overline{f_i}\big/\overline{h}
	\end{equation}
for $i=1,\dots,s$, where $\overline{f_i} = \frac{1}{n} \sum_{j=1}^n f_i({x}_j)$ and $\overline{h}\coloneqq \frac{1}{n}\sum_{j=1}^n h({x}_j)$.
	\end{itemize}
\end{theorem}
\vspace{0.2cm}

\begin{IEEEproof}
   See Sec-B.1 of Appendix.
\end{IEEEproof}
\vspace{0.2cm}

\begin{remark}
It is easy to see from the above theorem that $\overline{f}/\overline{h}\coloneqq \big[\overline{f_1}/\overline{h}, \dots, \overline{f_s}/\overline{h}\big]^\top$ is a sufficient statistic for the $k$-parameter $\mathbb{M}^{(\alpha)}$-family with respect to Jones et al. likelihood function.  
\end{remark}
\vspace{0.2cm}

\begin{example}
	\label{2expl:suff_stat_student_jones}
Consider the {\em Student distributions} with location parameter $\mu$, scale parameter $\sigma$ and degrees of freedom parameter $\nu\in(2,\infty):$
\begin{equation}
\label{1eqn:student_distribution_in_nu}
p_{\mu,\sigma} (x) = N_{\sigma,\nu} \bigg[ 1 +  \frac{(x-\mu)^2}{\nu\sigma^2} \bigg] ^{-\frac{\nu + 1}{2}},
\end{equation}	
for $x\in\mathbb{R}$, where $N_{\sigma,\nu} \coloneqq \dfrac{\Gamma(\frac{\nu+1}{2})}{\sqrt{\nu\pi\sigma^2}~\Gamma(\frac{\nu}{2})}$ is the normalizing factor. Observe that \eqref{1eqn:student_distribution_in_nu} includes \eqref{2eqn:student_distribution_alpha2} for $\nu=3$.

Let $\nu$ be known, and $\mu$ and $\sigma^2$ be unknown.
Letting $\alpha= \dfrac{\nu-1}{\nu+1}$ and $\theta = [\mu,\sigma^2]^\top$,
(\ref{1eqn:student_distribution_in_nu}) can be re-written as
\begin{eqnarray}
    \label{1eqn:student_t_as_B_alpha}
{p_{\theta} (x)
= N_{\theta,\alpha} \Big(1+\frac{\mu^2 b_\alpha}{\sigma^2}\Big)^{\frac{1}{\alpha -1}} }  \bigg[ 1 + \frac{b_\alpha x^2}{\sigma^2+b_\alpha\mu^2} -\frac{2\mu b_\alpha x}{\sigma^2+b_\alpha \mu^2}
\bigg]^{\frac{1}{\alpha -1}},~~
\end{eqnarray}
where $b_\alpha = \dfrac{1}{\nu} =\dfrac{1-\alpha}{1+\alpha}$ and 
$N_{\theta,\alpha}=N_{\sigma,\nu}$.
\vspace{0.2cm}

(\ref{1eqn:student_t_as_B_alpha}) implies that Student distributions form an  $\mathbb{M}^{(\alpha)}$-family with
\begin{align*}
\begin{array}{ll}
	h(x)\equiv 1,\quad Z(\theta) = N_{\theta,\alpha} \Big(1+\frac{\mu^2 b_\alpha}{\sigma^2}\Big)^{\frac{1}{\alpha -1}},\\
{w(\theta)=\big[w_1(\theta),w_2(\theta)\big]^\top,\quad
	f(x) = \big[f_1(x),
	f_2(x))\big]^\top,}\\
{w_1(\theta) = {b_\alpha }\big/{(\sigma^2+b_\alpha\mu^2)},\quad
	f_1(x) = x^2,} \\
{w_2(\theta) = {-2\mu b_\alpha }\big/{(\sigma^2+b_\alpha \mu^2)},\quad\text{and}\quad
	f_2(x) = x}.
\end{array}
\end{align*}
By Theorem \ref{2thm:suffiicent_statistics_M_alpha}, $\bar{f_1}= \frac{1}{n}\sum_{i=1}^n x_i^2$ and $\bar{f_2}= \frac{1}{n}\sum_{i=1}^n x_i$ are sufficient statistics for $\mu$ and $\sigma$ with respect to Jones et al. likelihood function.	
\end{example}


\subsection{Generalized Fisher-Darmois-Koopman-Pitman theorem for Basu et al. likelihood function}
\vspace{0.2cm}

In \cite{GayenK21ISIT}, we characterized the probability distributions that have a fixed number of sufficient statistics with respect to the Basu et al. likelihood function (only for the case $\Theta\subset \mathbb{R}$ though). We termed the associated family of probability distributions as $\mathbb{B}^{(\alpha)}$-family. In the following, we first define the $\mathbb{B}^{(\alpha)}$-family and then state the results concerning sufficient statistics for the same.
\vspace{0.2cm}

\begin{definition}
	\label{2defn:b_alpha_family}
	Let $h,w,f,\Theta$ and $\mathbb{S}$ be as in Definition \ref{2defn:general_M_alpha_family}. A {\em $k$-parameter $\mathbb{B}^{(\alpha)}$-family} characterized by $h,w,f,\Theta$ and $\mathbb{S}$ is defined by $\mathbb{B}^{(\alpha)}=\{p_\theta\colon \theta\in\Theta\}$, where
	\begin{eqnarray}
	\label{2eqn:form_of_general_b_alpha_family}
	{p_\theta(x)} = \left\{
	\begin{array}{ll}
	{[h(x) + F(\theta) +w(\theta)^\top f(x)]^{\frac{1}{\alpha -1}}} \text{~for~} 
	x \in \mathbb{S}\\
	{0} \hspace{4.5cm} \text{~otherwise}
	\end{array}
	\right.
	\end{eqnarray}
	for some  $F:\Theta\to\mathbb{R}$, the normalizing factor.
\end{definition}
\vspace{0.2cm}

\begin{remark}
\label{2rem:remark_on_b_alpha_family}
\begin{itemize}
    \item [(i)] Student distributions fall in $\mathbb{B}^{(\alpha)}$-family \cite[Ex. 1]{GayenK21J}.
    \vspace{0.2cm}
    
    \item[(ii)] $\mathbb{B}^{(\alpha)}$-family is also a generalization of the usual exponential family \cite{GayenK21J, GayenK21ISIT}. A more general form of this family can be found in \cite{CsiszarM12J} in connection to the projection theorem of Bregman divergences (c.f. \cite{GayenK21J}). This family was also studied in \cite{Bashkirov04J, Naudts04J}. A comparison between $\mathbb{B}^{(\alpha)}$ and $\mathbb{M}^{(\alpha)}$ families can be found in \cite[Rem. 2]{GayenK21J}. 
\end{itemize}
\end{remark}
\vspace{0.2cm}

The following is the generalized Fisher-Darmois-Koopman-Pitman theorem with respect to Basu et al. likelihood function $L_B^{(\alpha)}$.
\vspace{0.2cm}

\begin{theorem}
	\label{2thm:necessary_condition_sufficiency_basu}
	Let $k\leq s<n$. If $[T_1,\dots,T_s]^\top$ is a sufficient statistic for $\theta$ with respect to the Basu et al. likelihood function $L_B^{(\alpha)}$, then for $x\in \mathbb{S}$, there exist functions $h:\mathbb{S}\to  \mathbb{R}$, $f_i:\mathbb{S}\to \mathbb{R}$ and $w_i:\Theta\to \mathbb{R}$,  $i=1,\dots,s$ such that $p_\theta(x)$ can be expressed as (\ref{2eqn:form_of_general_b_alpha_family}).  
\end{theorem}
\vspace{0.2cm}

\begin{IEEEproof}
    The proof for the case $s=k$ can be found in \cite[Th. 2]{GayenK21ISIT}. On the other hand, let $k<s$. We give the sketch of the proof for $k=1$ and $s=2$ as in Theorem \ref{2thm:necessary_condition_sufficiency_jones}. The other cases can be proved similarly. Following the same steps as in Theorem \ref{2thm:necessary_condition_sufficiency_jones}, one gets (\ref{2eqn:likelihood_two_suff_stat_jones}) using $L_B^{(\alpha)}(x_1^n;\theta)$ in place of $L_J^{(\alpha)}(x_1^n;\theta)$. This, upon taking partial derivative with respect to $x_j$ first, and then with respect to some $x_i$, $i\neq j$, implies (\ref{2eqn:derivative}). Then, proceeding as in Theorem \ref{2thm:necessary_condition_sufficiency_jones}, one can prove the result.
\end{IEEEproof}
\vspace{0.2cm}

\noindent
The converse of the above theorem is also true. 
\vspace{0.2cm}

\begin{theorem}
	\label{2thm:suffiicent_statistics_b_alpha}
	Consider a $k$-parameter $\mathbb{B}^{(\alpha)}$-family as in Definition \ref{2defn:b_alpha_family}. Let $T=[T_1,\dots,T_s]^\top$ be $s$ statistics defined on $\mathbb{S}^n$ where $s\geq k$. Then $T$ is sufficient statistic for $\theta$ with respect to $L_B^{(\alpha)}$ if there exist real valued functions $\psi_1,\dots,\psi_s$ defined on $\mathscr{T}$ such that $\psi _i (T({x}_1^n)) = \overline{f_i}$ for any i.i.d. sample $x_1^n$ and $i=1,\dots,s$.
In particular, $\overline{f}\coloneqq \big[\overline{f_1}, \dots, \overline{f_s}\big]^\top$ is a sufficient statistic for $\theta$ with respect to $L_B^{(\alpha)}$.
\end{theorem}
\vspace{0.2cm}

\begin{IEEEproof}
   See \cite[Th. 4]{GayenK21ISIT}.
\end{IEEEproof}
\vspace{0.2cm}

\begin{example}
Consider the Student distributions as in Example \ref{2expl:suff_stat_student_jones}.
In \cite{GayenK21J}, we showed that these also form a 2-parameter $\mathbb{B}^{(\alpha)}$-family. Using Theorem \ref{2thm:suffiicent_statistics_b_alpha}, the same statistics as in Example \ref{2expl:suff_stat_student_jones} form sufficient statistics for their parameters with respect to Basu et al. likelihood function.
\end{example}

\subsection{Minimal sufficiency with respect to Basu et al. and Jones et al. likelihood functions}	\vspace{0.2cm}

In the following, we show that $\bar{f}/\bar{h}$ and $\bar{f}$ are also minimal, respectively, for regular $\mathbb{M}^{(\alpha)}$ and $\mathbb{B}^{(\alpha)}$ families.
An $\mathbb{M}^{(\alpha)}$-family as in Definition \ref{2defn:general_M_alpha_family} is said to be \emph{regular} if the following are satisfied (c.f. \cite[Def.~3]{GayenK21J}).
\begin{itemize}
    \item[(i)] The support $\mathbb{S}$ does not depend on the parameter $\theta$.
    \item[(ii)] The number of unknown parameters is the same as the number of $w_i$'s, that is, $s=k$.
    \item[(iii)] $1,w_1,\dots, w_k$ are linearly independent.
    \item[(iv)] $f_1,\dots, f_k$ are linearly independent.
\end{itemize}
\vspace{0.2cm}

Similarly, a $\mathbb{B}^{(\alpha)}$-family as in Definition \ref{2defn:b_alpha_family} is said to be \emph{regular} if $1, f_1,\dots, f_k$ are linearly independent along with conditions (i)-(iii) above (c.f. \cite[Def. 1]{GayenK21J}).
\vspace{0.2cm}

\begin{theorem}
	\label{2thm:minimal_sufficient_statistics_M_alpha}
	The statistics $\overline{f}/\overline{h} = \big[\overline{f_1}/\overline{h}, \dots, \overline{f_k}/\overline{h}\big]^\top$ and $\overline{f} = \big[\overline{f_1}, \dots, \overline{f_k}\big]^\top$ are 
	minimal sufficient statistics, respectively, for a $k$-parameter regular $\mathbb{M}^{(\alpha)}$-family and a $k$-parameter regular $\mathbb{B}^{(\alpha)}$-family with respect to Jones et al. and Basu et al. likelihood functions.
\end{theorem}
\vspace{0.2cm}

\begin{IEEEproof}
   See Sec-B.2 of Appendix.
\end{IEEEproof}
\vspace{0.2cm}

\begin{remark}
	In view of Theorem \ref{2thm:minimal_sufficient_statistics_M_alpha}, the sufficient statistics derived in Example \ref{2expl:suff_stat_student_jones} for Student distributions are minimal too (with respect to both Basu et al. and Jones et al. likelihood functions).
\end{remark}
	\section{Rao-Blackwell theorem for Jones et al. estimation}
\label{2sec:Rao_blackwelLones_et_al}

Let $X_1^n=x_1^n$ be an i.i.d. sample drawn according to some $p_{\theta}\in \Pi$. Here and afterward, we assume  $\theta\in\Theta\subset \mathbb{R}$ unless otherwise stated. Let $\widehat{\theta}\coloneqq\widehat{\theta} (x_1^n)$ be an unbiased estimator of $\theta$. Let $T\coloneqq T(x_1^n)$ be a sufficient statistic for $\theta$ with respect to the usual log-likelihood function. 
Let us consider the conditional expectation $\mathbb{E}_\theta[\widehat{\theta}|T]$. It is easy to show from Fisher's definition of sufficiency that $\mathbb{E}_\theta[\widehat{\theta}|T]$ is independent of $\theta$ \cite{Blackwell47J}.
Hence $\mathbb{E}_\theta[\widehat{\theta}|T]$ is  a function of $T$ alone. Let $\phi(T)\coloneqq \mathbb{E}_\theta[\widehat{\theta}|T]$. Then by Rao-Blackwell theorem, $\phi(T)$ is a {\em uniformly better} unbiased estimator of $\theta$ than $\widehat{\theta}$ (see \cite{Rao45J, Blackwell47J}). In other words,
\begin{itemize}

    \item [1.] $\mathbb{E}_\theta[\phi(T)] = \mathbb{E}_\theta[\widehat{\theta}] = \theta$, and
    
    \item [2.] $\text{Var}_\theta[\phi(T)]\leq \text{Var}_\theta[\widehat{\theta}]$ 
    
\end{itemize}
for all $\theta\in\Theta$, and equality holds if and only if $\widehat{\theta}$ is a function of the sufficient statistic $T$.
This implies that the `best' unbiased estimator of $\theta$ (the one with the least variance) is a function of the sufficient statistic $T$ only, that is, one from the family $\{\psi(T): \mathbb{E}_\theta[\psi(T)] = \theta\}$. 
%

Let us consider an exponential family $\mathcal{E}$ as in (\ref{2eqn:form_of_exp_family}).
$\bar{f} = \frac{1}{n}\sum_{i=1}^n f(x_i)$ is a sufficient statistic for $\theta$ with respect to the usual log-likelihood function \cite{LehmannC98B} and $\mathbb{E}_{\theta}[\bar{f}] =\mathbb{E}_{\theta}[f(X)]=: \tau(\theta)$, called the \emph{mean value parameter}. 
Let $\widehat{\theta}$ be any unbiased estimator of $\tau(\theta)$. Then by Rao-Blackwell theorem, we have $\phi(\bar{f}) \coloneqq\mathbb{E}_\theta[\widehat{\theta}|\bar{f}]$ is also an unbiased estimator of $\tau(\theta)$ and has variance not more than that of $\widehat{\theta}$. That is,
\begin{eqnarray}
    \label{2eqn:rao_blackwell_exp_model}
\mathbb{E}_\theta[\phi(\bar{f})] = \tau(\theta)\quad \text{and}\quad\text{Var}_\theta[\phi(\bar{f})] \leq \text{Var}_\theta[\widehat{\theta}].
\end{eqnarray}
Now consider $\psi(\bar{f})=\phi(\bar{f}) - \bar{f}$. Then 
\begin{equation*}
   \mathbb{E}_\theta[\psi(\bar{f})] = \mathbb{E}_\theta[\phi(\bar{f})] - \mathbb{E}_\theta[\bar{f}] = \tau(\theta) - \tau(\theta) =0, 
\end{equation*}
and
\begin{eqnarray}
    \label{var_inq}
{\text{Var}_\theta[\phi(\bar{f})]}
&=& \text{Var}_\theta[\psi(\bar{f}) + \bar{f}]\nonumber\\
&=& \text{Var}_\theta[\psi(\bar{f})] + \text{Var}_\theta[\bar{f}] + 2~ \text{Cov}_\theta[\psi(\bar{f}),\bar{f}].~~~~
\end{eqnarray}
 It can be shown, under some mild regularity conditions, that $\text{Cov}_\theta[\psi(\bar{f}),\bar{f}] = 0$. See Sec-C.1 of Appendix for the detailed derivation. Hence (\ref{var_inq}) reduces to 
\begin{equation*}
\text{Var}_\theta[\phi(\bar{f})] =  \text{Var}_\theta[\psi(\bar{f})] + \text{Var}_\theta[\bar{f}] \geq \text{Var}_\theta[\bar{f}].
\end{equation*}
Thus from (\ref{2eqn:rao_blackwell_exp_model}), for any $\widehat{\theta}$ with $\mathbb{E}_\theta[\widehat{\theta}] = \mathbb{E}_\theta[\bar{f}]$, we have
\begin{equation*}
\text{Var}_\theta[\widehat{\theta}]\geq \text{Var}_\theta[\phi(\bar{f})]
\geq \text{Var}_\theta[\bar{f}].
\end{equation*}
This shows that the variance of $\bar{f}$ is not more than that of any unbiased estimator of $\tau(\theta)$. Thus we have the following (c.f. \cite{Rao49J, CasellaB02B}).
\vspace{0.2cm}

\begin{proposition}
\label{2prop:rao_blackwell_best_est}
	Consider a regular exponential family as in (\ref{2eqn:form_of_exp_family}). Let $\mathbb{E}_\theta[f(X)] = \tau(\theta)$. Then there does not exist any unbiased estimator of $\tau(\theta)$ that has lesser variance than that of $\bar{f}$. That is, $\bar{f}$ is the best estimator of its expected value.
\end{proposition} 
\vspace{0.2cm}

\noindent
Further, it can be easily shown that (c.f. \cite{Rao47J, Rao49J})
\begin{equation*}
    \text{Var}_{\theta}[\bar{f}]= \frac{1}{n}\text{Var}[f(X)] = \frac{[\tau'(\theta)]^2}{n~I(\theta)},
\end{equation*}
where $I(\theta)$ is the Fisher information, that is, $I(\theta)=\mathbb{E}_\theta\Big[\big(\frac{\partial}{\partial\theta}\log p_\theta(X)\big)^2\Big]$. This and Proposition \ref{2prop:rao_blackwell_best_est} imply that
	\begin{equation}
	\label{2eqn:CR_inq_exp_model}
	\text{Var}_{\theta}[\widehat{\theta}]\geq \frac{[\tau'(\theta)]^2}{n~I(\theta)}.
\end{equation}
Equality holds in \eqref{2eqn:CR_inq_exp_model} if and only if $\widehat{\theta}=\bar{f}$ almost surely. Notice that (\ref{2eqn:CR_inq_exp_model}) is the well-known \emph{Cram\'er-Rao inequality} \cite{Rao45J, Cramer46B}. 

 Observe that, in Rao-Blackwell theorem, instead of a sufficient statistic $T$, if one conditions on some other statistic, then also the resulting estimator has variance not more than that of the given estimator. However, conditioning on any other statistic does not guarantee that the conditional expectation is independent of $\theta$. Also by Fisher-Darmois-Koopman-Pitman theorem, a sufficient statistic with respect to the usual log-likelihood function exists if and only if the underlying model is exponential \cite{Koopman36J, Pitman36J, Darmois35J}. Hence Rao-Blackwell theorem is not applicable outside the exponential family. However, in Section \ref{2sec:suff_stat_jones}, we showed that there are models that have a fixed number of sufficient statistics with respect to likelihood functions other than log-likelihood. In particular, we showed that $\mathbb{M}^{(\alpha)}$-family is one such model which has a fixed number of sufficient statistics with respect to Jones et al. likelihood function (see Theorem \ref{2thm:necessary_condition_sufficiency_jones} and Theorem \ref{2thm:suffiicent_statistics_M_alpha}). In the following, we generalize the Rao-Blackwell theorem for Jones et al. likelihood function. Later we will apply this to $\mathbb{M}^{(\alpha)}$-family and find Rao-Blackwell type ``best'' estimator. This also helps us find a lower bound for the variance of the estimators of $\mathbb{M}^{(\alpha)}$-family as an alternative to the usual CRLB.

The Jones et al. likelihood function $L_J^{(\alpha)}$ as in (\ref{2eqn:likelihood_function_for_I_alpha}) can be re-written as
\begin{eqnarray*}
    L_J^{(\alpha)}(x_1^n ; \theta)
    &=& \log \Big[{\Big(\frac{1}{n} \sum\limits_{i=1}^n p_\theta(x_i)^{\alpha -1}\Big)^{\frac{1}{\alpha -1}}}/{\|p_\theta\|}\Big].
\end{eqnarray*}
Then the deformed probability distribution associated with $L_J^{(\alpha)}$ is given by
\begin{eqnarray}
\label{gen_joint_distribution_jones_et_al}
    p_\theta^*(x_1^n) &\coloneqq& \frac{\exp[ L_J^{(\alpha)}(x_1^n; \theta)]}{\int \exp[ L_J^{(\alpha)}(y_1^n; \theta)] dy_1^n}\nonumber\\
    &~=& \dfrac{\bigg(\frac{1}{n} \sum\limits_{i=1}^n p_\theta(x_i)^{\alpha -1}\bigg)^{\frac{1}{\alpha -1}}}{\int \bigg(\frac{1}{n} \sum\limits_{i=1}^n p_\theta(y_i)^{\alpha -1}\bigg)^{\frac{1}{\alpha -1}} dy_1^n}.
\end{eqnarray}
Observe that when $\alpha\to 1$, $L_J^{(\alpha)}(x_1^n;\theta)$ coincides with the log-likelihood function $L(x_1^n;\theta)$. Hence, in this case,  $p_\theta^*(x_1^n) =  p_\theta(x_1^n)$.

Let $T$ be a statistic of the sample $x_1^n$. Then, from (\ref{2eqn:gen_condition_dist}), we have
 \[ p_{\theta{_{X_1^n|T}}}^*(x_1^n|t) =
 \begin{cases} 
      \dfrac{p_\theta^*(x_1^n)}{q_\theta^*(t)} & \text{if}~x_1^n\in A_t \\
      0 & \text{otherwise}. 
   \end{cases}
\]
Observe that, as $\alpha\to 1$, the above coincides with
the usual conditional distribution of $x_1^n$ given $T=t$.

Let $\widehat{\theta}\coloneqq\widehat{\theta}(x_1^n)$ be an unbiased estimator of $\tau^*(\theta)\coloneqq\mathbb{E}_\theta^*[\widehat{\theta}]$, where $\mathbb{E}_\theta^*[\cdots]$ denotes the expectation with respect to $p_\theta^*$. Let $T\coloneqq T(x_1^n)$ be a sufficient statistic for $\theta$ with respect to Jones et al. likelihood function $L_J^{(\alpha)}$.
 Let us consider the conditional expectation
\begin{eqnarray}
    \label{2eqn:gen_cond_distr}
    \mathbb{E}_\theta^*\big[\widehat{\theta}|T=T(x_1^n)\big]
    = \int \widehat{\theta}(y_1^n) p_\theta^*(y_1^n|T(x_1^n))dy_1^n.
\end{eqnarray}
In the following, we first establish some properties of  (\ref{2eqn:gen_cond_distr}).
\vspace{0.2cm}

\begin{proposition}
   \label{thm1}
	The conditional expectation $\mathbb{E}_{\theta}^*[\widehat{\theta}|T]$ is independent of $\theta$. 
\end{proposition}
\vspace{0.2cm}

\begin{IEEEproof}
	See Sec-C.2 of Appendix.
\end{IEEEproof}
\vspace{0.2cm}

\noindent
Let
\begin{equation*}
    \mathbb{E}_\theta^*\big[\widehat{\theta}|T=T(x_1^n)\big] \eqqcolon \phi^*(T(x_1^n)).
\end{equation*}
Let us consider the expectation of $\phi^*(T)$ with respect to $p_\theta^*$. That is,
\begin{equation*}
    \mathbb{E}_\theta^*[\phi^*(T(X_1^n))] = \int p_\theta^*(x_1^n) \phi^*(T(x_1^n))dx_1^n.
\end{equation*}
Then we have the following.
\vspace{0.2cm}

\begin{proposition}
   \label{thm2}
	$\phi^*(T)$ is unbiased for $\tau^*(\theta)$, that is, $\mathbb{E}_\theta^*[\phi^*(T(X_1^n))] =   \mathbb{E}_\theta^*[\widehat{\theta}(X_1^n))]=\tau^*(\theta).$
\end{proposition}
\vspace{0.2cm}

\begin{IEEEproof}
	See Sec-C.3 of Appendix.
\end{IEEEproof}
\vspace{0.2cm}

\noindent
Thus, if we have an unbiased estimator of $\tau^*(\theta)$ and a sufficient statistic with respect to Jones et al. likelihood function, we can always get another unbiased estimator for $\tau^*(\theta)$ that is based only on the sufficient statistic. Let us now consider the following two classes of estimators based on the sufficient statistic $T:$
\begin{eqnarray}
\label{2eqn:unbiased_estimator_based_on_suffi_stat}
\mathcal{A} &=& \{\mathbb{E}^*_{\theta}[\widehat{\theta}|T]:
	\quad\widehat{\theta}~\text{is an unbiased estimator of}\tau^*(\theta),~\text{that is,}~\mathbb{E}_\theta^*[\widehat{\theta}] = \tau^*(\theta)\},\nonumber\\
 \mathcal{B} &=& \{\psi(T): \mathbb{E}^*_{\theta}[\psi(T)] = \tau^*(\theta)\}.
\end{eqnarray}

\begin{proposition}
   \label{thm3}
	The above two classes of estimators are the same, that is $\mathcal{A}=\mathcal{B}$.
\end{proposition}
\vspace{0.2cm}

\begin{IEEEproof}
    See Sec-C.4 of Appendix.
\end{IEEEproof}
\vspace{0.2cm}

The above proposition suggests that any unbiased estimator of $\tau^*(\theta)$ that is a function of the sufficient statistic $T$ is of the form $\mathbb{E}^*_{\theta}[\widehat{\theta}|T]$, where $\widehat{\theta}$ is an unbiased estimator of $\tau^*(\theta)$. In the following, we compare the variance of an estimator of the form (\ref{2eqn:unbiased_estimator_based_on_suffi_stat}) and any other unbiased estimator of $\tau^*(\theta)$.
\vspace{0.2cm}

\begin{lemma}
    \label{thm4}
	Let $\widehat{\theta}$ be any unbiased estimator of $\tau^*(\theta)$ and $\phi^*(T) =\mathbb{E}_\theta^*[\widehat{\theta}|T]$. Let $\psi^*(T)\in\mathcal{B}$. Then
	\begin{eqnarray}
	    \label{variance_of_any_esti}
		{\text{Var}^*_{\theta}[\widehat{\theta}] = \mathbb{E}^*_{\theta}[(\widehat{\theta} - \psi^*(T))^2] + \text{Var}^*_{\theta}[\psi^*(T)]} + 2 \mathbb{E}^*_{\theta}[\psi^*(T)\{\phi^*(T) -\psi^*(T)\}].~~
	\end{eqnarray}
\end{lemma}
\vspace{0.2cm}

	\begin{IEEEproof}
	   See Sec-C.5 of Appendix.
	\end{IEEEproof}
\vspace{0.2cm}

The following (also the main result of this section) is a generalization of Rao-Blackwell theorem for Jones et al. likelihood function. 
\vspace{0.2cm}

\begin{theorem}
    \label{2thm:gen_rao_blackwell}
Let $\Pi=\{p_\theta:\theta\in\Theta\subset \mathbb{R}\}$ be a family of probability distributions with common support $\mathbb{S}\subset\mathbb{R}$. Consider the deformed family of probability distributions $\Pi^*=\{p_\theta^*\}$ of $\Pi$, where $p_\theta^*$ is as in (\ref{gen_joint_distribution_jones_et_al}).
Let $\widehat{\theta}$ be an unbiased estimator for its expected value with respect to $p_\theta^*$. Let $\tau^*(\theta) \coloneqq \mathbb{E}_\theta^*[\widehat{\theta}]$. Suppose $T$ is a sufficient statistic for $\theta$ with respect to the Jones et al. likelihood function. Let $\phi^*(T) \coloneqq\mathbb{E}_\theta^*[\widehat{\theta}|T]$. Then $\phi^*(T)$ is a uniformly better estimator of $\tau^*(\theta)$ than $\widehat{\theta}$, that is,

\begin{itemize}
    \item [(a)] $\mathbb{E}_\theta^* [\phi^*(T)] = \tau^*(\theta)$,
    
    \item[(b)] $\text{Var}_\theta^*[\phi^*(T(X_1^n))]\leq \text{Var}_\theta^*[\widehat{\theta}(X_1^n)]$, and equality holds if and only if $\widehat{\theta}$ is a function of the sufficient statistic $T$.
\end{itemize}
\end{theorem}
\vspace{0.2cm}

\begin{IEEEproof}
	See Sec-C.6 of Appendix.
\end{IEEEproof} 
\vspace{0.2cm}

Theorem \ref{2thm:gen_rao_blackwell} implies that $\phi^*(T)$ is \emph{uniformly a better estimator}  of $\tau^*(\theta)$ than $\widehat{\theta}$. This means that the \emph{best unbiased estimator} of $\tau^*(\theta)$ is a function of the sufficient statistic $T$.
The following result asserts that the estimator with minimum variance (best estimator) is unique.
\vspace{0.2cm}

\begin{proposition}
\label{unique_best_est_wrt_p_star}
The best unbiased estimator for $\tau^*(\theta)$ with respect to $p_\theta^*$ is unique.
\end{proposition}
\vspace{0.2cm}

\begin{IEEEproof}
    See Sec-C.7 of Appendix.
\end{IEEEproof}

In the following, we apply the above results to the family of Bernoulli distributions.
\subsection{Example:}
Observe that the family of Bernoulli distributions can be written as 
\begin{equation*}
    p_\theta(x) = (1-\theta) [1+x\cdot (2\theta-1)/(1-\theta)]\quad\text{for}\quad x=0,1,
\end{equation*}
where $\theta\in(0,1)$ is the unknown parameter.
Thus Bernoulli distributions form an $\mathbb{M}^{(\alpha)}$-family with $\alpha =2$, $h(x)\equiv 1$, $Z(\theta)=1-\theta$, $w(\theta)=(2\theta -1)/(1-\theta)$ and $f(x)=x$. By Theorem \ref{2thm:suffiicent_statistics_M_alpha}, $\sum_{i=1}^n x_i$ is a sufficient statistic for $\theta$ with respect to the Jones et al. likelihood function.

The Jones et al. likelihood function for Bernoulli distributions is
\begin{eqnarray*}
L_{J}^{(2)}(x_1^n,\theta)
&=& \log \Bigg[\frac{\frac{1}{n}\sum_{i=1}^n p_\theta(x_i)}{(\mathbb{E}_\theta[p_\theta(X)])^{1/2}}\Bigg]\\
&=& \log\Bigg[ \frac{(1-\theta)^{1/2}\{1+\overline{x}\cdot (2\theta-1)/(1-\theta)\} }{(1+\theta(2\theta-1)/(1-\theta))^{1/2}}\Bigg].
\end{eqnarray*}
The associated deformed probability distribution is given by
\begin{eqnarray*}
p_\theta^*(x_1^n) &=& \dfrac{\exp[L^{(2)}_{J}(x_1^n,\theta)]}{\sum_{y_1^n}\exp[L^{(2)}_{J}(y_1^n,\theta)]}\\
&=&\frac{1+ \overline{x}\cdot (2\theta-1)/(1-\theta)}{\sum_{y_1^n}[1+ \overline{y}\cdot (2\theta-1)/(1-\theta)]}\\
&=& N(\theta) \bigg[1+\overline{x}\cdot \frac{2\theta-1}{1-\theta}\bigg],
\end{eqnarray*}
where 
\begin{eqnarray*}
N(\theta)^{-1}&\coloneqq& \sum_{y_1^n}\bigg[1+\overline{y}\cdot \frac{2\theta-1}{1-\theta}\bigg]\nonumber\\
&=& 2^n + \left[\frac{(2\theta-1)}{(1-\theta)}\cdot \frac{1}{n}\sum_{i=1}^n i\cdot\binom{n}{i}\right]\nonumber\\
&=& 2^n + \bigg[\frac{(2\theta-1)}{n(1-\theta)} \cdot \big(n2^{n-1}\big)\bigg]\nonumber\\
&=& \frac{2^{n-1}}{1-\theta}.
\end{eqnarray*}
Also
\begin{eqnarray*}
    \lefteqn{ p_\theta^*\Bigg(X_1^n=x_1^n\Big|T=\sum_{i=1}^n x_i\Bigg)} \\
     &=& \dfrac{1+\overline{x}\cdot {(2\theta-1)}/{(1-\theta)}}{\sum\limits_{y_1^n\in A_{\sum x_i}}[1+\overline{y}\cdot (2\theta-1)/(1-\theta)]}\\
&=& \dfrac{1}{\big|A_{\sum x_i}\big|}\\
&=& \frac{1}{\binom{n}{\sum x_i}}.
\end{eqnarray*}
Let $\widehat{\theta} = x_1$. Then
\begin{eqnarray*}
\lefteqn{\tau^*(\theta)}\\
&=& \mathbb{E}_\theta^*[\widehat{\theta}]\\
&=& \sum\limits_{x_1^n} x_1 \cdot p_\theta^*(x_1^n)\\
&=& \sum\limits_{x_1} N(\theta)\Bigg[2^{n-1} +\frac{2\theta-1}{n(1-\theta)}\cdot \sum_{i=1}^n (x_1+i)\binom{n-1}{i}\Bigg]\nonumber\\
&=& N(\theta)\cdot 2^{n-1} \bigg[1 + \frac{2\theta-1}{n(1-\theta)}\cdot \bigg(1+\frac{n+1}{2}\bigg)\bigg]\nonumber\\
&=& \frac{\theta}{n} + \frac{n-1}{2n}.
\end{eqnarray*}
The variance of $\widehat{\theta}$ with respect to $p_\theta^*$ is given by
\begin{eqnarray*}
\text{Var}_\theta^*[\widehat{\theta}] &=& \mathbb{E}_\theta^*\big[X_1^2\big] - \big(\mathbb{E}_\theta^*[X_1]\big)^2
= \tau^*(\theta) - \tau^*(\theta)^2.
\end{eqnarray*}
Let us consider the generalized Rao-Blackwell estimator of $\tau^*(\theta)$. Recall that $\sum_{i=1}^n x_i$ is a sufficient statistic for $\theta$ with respect to Jones et al. likelihood function.
\begin{eqnarray*}
\phi^*\big(T(x_1^n)\big)&=&\phi^*\left(\sum_{i=1}^n x_i\right)\nonumber\\
&=& \sum\limits_{y_1^n} \widehat{\theta}(y_1^n)\cdot~ p_\theta^*\left(y_1^n~\Big|~T=\sum_{i=1}^n x_i\right)\nonumber\\
&=& \sum\limits_{y_1^n\in A_{\sum x_i}} y_1 \cdot~ p_\theta^*\Bigg(y_1^n~\Big|~T=\sum_{i=1}^n x_i\Bigg)\nonumber\\
&=& {\binom{n}{\sum x_i}}^{-1} \cdot \sum\limits_{y_1^n\in A_{\sum x_i}} y_1\nonumber\\
&=& {\binom{n}{\sum x_i}}^{-1} \cdot \sum\limits_{y_2^n\in A_{(\sum_{i=1}^n x_i -1)}} 1\nonumber\\
&=& \binom{n-1}{\sum x_i -1}\Bigg/ \binom{n}{\sum x_i}\nonumber\\
&=& \frac{1}{n}\sum_{i=1}^n x_i.
\end{eqnarray*}
Thus $\phi^*(T(x_1^n)) =\frac{1}{n}\sum\limits_{i=1}^n x_i$ is uniformly better estimator for $\tau^*(\theta)$ than $x_1$.
It can also be shown that
\begin{equation*}
    \text{Var}_\theta^*[\phi^*(T)] = \frac{\Big[\frac{\partial}{\partial\theta}\tau^*(\theta)\Big]^2}{I_{\alpha,n}^*(\theta)},
\end{equation*}
where $I_{\alpha,n}^*(\theta)$ is the inverse of the asymptotic variance of Jones et al. estimator for $p_\theta^*(X_1^n)$ \cite{JonesHHB01J}. This shows that the Rao-Blackwell estimator derived from a sufficient statistic attains the asymptotic variance of Jones et al. estimator. In the following section, we elucidate this for a general $\mathbb{M}^{(\alpha)}$-family. 
\section{Efficiency of Rao-Blackwell estimators of power-law family}
\label{2sec:application_student_dist}
In this section, we shall show the efficiency of Rao-Blackwell estimators for a power-law family. For simplicity, we shall consider a regular $\mathbb{M}^{(\alpha)}$-family with $h\equiv 1$. That is, for $x\in \mathbb{S}$,
\begin{equation}
\label{2eqn:m_alpha_family}
    p_\theta(x) = Z(\theta) [1+w(\theta)f(x)]^{\frac{1}{\alpha -1}}.
\end{equation}
Then
\begin{eqnarray*}
   {\left(\frac{1}{n}\sum\limits_{i=1}^n p_\theta(x_i)^{\alpha -1} \right)^{\frac{1}{\alpha -1}}}
    &=& \left(\frac{1}{n}\sum\limits_{i=1}^n Z(\theta)^{\alpha -1} [1 + w(\theta) f(x_i)]\right)^{\frac{1}{\alpha -1}}\nonumber\\
    &=& Z(\theta) [1 +w(\theta)\overline{f}]^{\frac{1}{\alpha -1}}.
\end{eqnarray*}
Thus, using (\ref{gen_joint_distribution_jones_et_al}), we get
\begin{eqnarray}
\label{rao_blackwell_distri_jones_et}
    p_\theta^*(x_1^n) &=& \frac{[1 +w(\theta)\overline{f}]^{\frac{1}{\alpha -1}}}{\int \big[1+w(\theta) \bar{f}\big]^{\frac{1}{\alpha -1}}dx_1^n}
    = N(\theta) [1+w(\theta)\overline{f}]^{\frac{1}{\alpha -1}},
\end{eqnarray}
where $\bar{f}\coloneqq\frac{1}{n}\sum_{i=1}^n f(x_i)$ and $N(\theta)^{-1}\coloneqq \int [1+w(\theta) \bar{f}]^{\frac{1}{\alpha -1}}dx_1^n$ is the normalizing constant. Comparing this with (\ref{2eqn:form_of_general_M_alpha_family}), we see that (\ref{rao_blackwell_distri_jones_et}) also forms an $\mathbb{M}^{(\alpha)}$-family with $h\equiv 1$ and $f(x_1^n) =\bar{f}$. Recall that, $\bar{f}$ is a sufficient statistic for $\mathbb{M}^{(\alpha)}$ with respect to Jones et al. likelihood function (Theorem \ref{2thm:suffiicent_statistics_M_alpha}). Let $\tau^*(\theta) \coloneqq\mathbb{E}_\theta^*[\bar{f}] $. In the following, we show that $\bar{f}$ is the best estimator of $\tau^*(\theta)$. That is the variance of any unbiased estimator of $\tau^*(\theta)$ is bounded below by the variance of $\bar{f}$ with respect to $p_\theta^*$. 

Let $\widehat{\theta}$ be an unbiased estimator of $\tau^*(\theta)$ with respect to $p_\theta^*$. Then by Theorem \ref{2thm:gen_rao_blackwell}, $\phi^*(\bar{f})=\mathbb{E}_\theta^*[\widehat{\theta}|\bar{f}]$ is also an unbiased estimator of $\tau^*(\theta)$ and
\begin{equation}
\label{gen_rao_blackwell_variance}
\text{Var}^*_\theta[\widehat{\theta}] \geq \text{Var}^*_\theta[\phi^*(\bar{f})].
\end{equation} 
Let us consider $\psi^*(\bar{f})\coloneqq \phi^*(\bar{f}) - \bar{f}$. Then $\mathbb{E}_\theta^*[\psi^*(\bar{f})] = \tau^*(\theta) - \tau^*(\theta) = 0.$ Also
\begin{eqnarray}
    \label{2eqn:var_of_phi*}
{\text{Var}_\theta^*[\phi^*(\bar{f})]}
&=& \text{Var}_\theta^*[\psi^*(\bar{f}) + \bar{f}]\nonumber\\
&=& \text{Var}_\theta^*[\psi^*(\bar{f})] + \text{Var}_\theta^*[\bar{f}] + 2~ \text{Cov}_\theta^*[\psi^*(\bar{f}), \bar{f}].~~~~~
\end{eqnarray}
In Sec-D.1 of Appendix, we have shown that, if $w(\theta)>0$ for $\theta\in\Theta$, then $\text{Cov}_\theta^*[\psi^*(\bar{f}), \bar{f}] \geq 0$. Using this in (\ref{2eqn:var_of_phi*}), we get
\begin{eqnarray}
\label{suffi_stat_lesser_var_generalized}
{\text{Var}_\theta^*[\phi^*(\bar{f})]}
&=& \text{Var}_\theta^*[\psi^*(\bar{f})] + \text{Var}_\theta^*[\bar{f}] + 2~ \text{Cov}_\theta^*[\psi^*(\bar{f}), \bar{f}]\nonumber\\
&\geq& \text{Var}_\theta^*[\bar{f}].
\end{eqnarray}
Using (\ref{gen_rao_blackwell_variance}) and (\ref{suffi_stat_lesser_var_generalized}), we get
\begin{equation}
\label{2eqn:gen_cr_ineq}
\text{Var}^*_\theta[\widehat{\theta}] \geq  \text{Var}_\theta^*[\bar{f}],
\end{equation}
where equality holds if and only if $\widehat{\theta}=\bar{f}$ almost surely. Thus we have the following.
\vspace{0.2cm}

\begin{theorem}
	\label{2thm:best_est_m_alpha}
	Consider an $\mathbb{M}^{(\alpha)}$-family as in (\ref{2eqn:m_alpha_family}) with $w(\theta)> 0$ for $\theta\in\Theta\subset\mathbb{R}$. Let $\tau^*(\theta) = \mathbb{E}^*_\theta[\bar{f}]$. Then the variance of any unbiased estimator of $\tau^*(\theta)$ with respect to $p_\theta^*$ is not lower than that of $\bar{f}$. That is, $\bar{f}$ is the best estimator of its expected value with respect to $p_\theta^*$.
\end{theorem}
\vspace{0.2cm}

We now show that $\text{Var}_\theta^*[\bar{f}]$, indeed, equals to the asymptotic variance of Jones et al. estimator for the model $p_\theta^*$. 
The Jones et al. asymptotic variance for $p_\theta^*$ is given by (c.f. \cite{JonesHHB01J})
\begin{eqnarray}
    \label{2eqn:gen_fisher_information}
    \frac{1}{I_{\alpha,n}^*(\theta)}
    \coloneqq \frac{\text{Var}_\theta^*[p_\theta^*(X_1^n)^{\alpha -1}s^*(X_1^n,\theta)]}{\text{Cov}_\theta^*[s^*(X_1^n,\theta), p_\theta^*(X_1^n)^{\alpha -1}s^*(X_1^n,\theta)]^2},
\end{eqnarray}
where
\begin{equation}
    \label{2eqn:gen_score}
    s^*(x_1^n,\theta) \coloneqq \frac{\partial}{\partial\theta}\log p_\theta^*(x_1^n).
\end{equation}
It can be easily seen that $I_{\alpha,n}^*(\theta)$ coincides with the usual Fisher information $I_n^*(\theta)$ for $\alpha = 1$, where $I_n^*(\theta)= \text{Var}_\theta^*[\frac{\partial}{\partial\theta}\log p_\theta^*(X_1^n)]$.
Observe that (\ref{rao_blackwell_distri_jones_et}) can be re-written as
\begin{eqnarray*}
p_\theta^*(x_1^n)
    &=& N(\theta) [1+w(\theta)\overline{f}]^{\frac{1}{\alpha -1}}
    = \big[N(\theta)^{\alpha -1} + N(\theta)^{\alpha -1}w(\theta) \overline{f}\big]^{\frac{1}{\alpha -1}}.
\end{eqnarray*}
For convenience, let $F(\theta) \coloneqq N(\theta)^{\alpha -1}$ and $\widetilde{w}(\theta) \coloneqq N(\theta)^{\alpha -1}w(\theta)$. Then
\begin{equation*}
    p_\theta^*(x_1^n)
    = [F(\theta) + \widetilde{w}(\theta) \overline{f}]^{\frac{1}{\alpha -1}}
\end{equation*}
and hence
\begin{equation*}
    \frac{\partial}{\partial\theta}p_\theta^*(x_1^n)
    = \frac{1}{\alpha -1} p_\theta^*(x_1^n)^{2-\alpha} [F'(\theta) + \widetilde{w}'(\theta) \overline{f}].
\end{equation*}
Then 
\begin{eqnarray}
    \label{2eqn:var}
{\text{Var}_\theta^*\big[p_\theta^*(X_1^n)^{\alpha -1}s^*(X_1^n,\theta)\big]}
&=& \left(\frac{1}{\alpha -1}\right)^2 \big[\widetilde{w}'(\theta)\big]^2 \text{Var}^*_\theta [\bar{f}],
\end{eqnarray}
and
\begin{eqnarray}
\label{2eqn:covar}
 \text{Cov}^*_\theta[s^*(X_1^n,\theta), p_\theta^*(X_1^n)^{\alpha -1}s^*(X_1^n,\theta)]& = &\left(\frac{1}{\alpha -1}\right)\widetilde{w}'(\theta) \frac{\partial}{\partial\theta}\tau^*(\theta).
\end{eqnarray}
The detailed derivations of (\ref{2eqn:var}) and (\ref{2eqn:covar}) can be found in Sec-D.2 of the Appendix. Using these values, (\ref{2eqn:gen_fisher_information}) becomes
\begin{align*}
\frac{1}{I^*_{\alpha,n}(\theta)} = \frac{\text{Var}^*_\theta \big[\overline{f}\big]}{\Big[\frac{\partial}{\partial\theta}\mathbb{E}^*_\theta[\overline{f}]\Big]^2},
\end{align*}
\vspace{-0.4cm}

that is,
\vspace{-0.4cm}

\begin{eqnarray*}
 \dfrac{{\Big[\frac{\partial}{\partial\theta}\mathbb{E}^*_\theta[\overline{f}]\Big]^2}}{I^*_{\alpha,n}(\theta)} =
 \frac{{\Big[\frac{\partial}{\partial\theta}\tau^*(\theta)\Big]^2}}{I^*_{\alpha,n}(\theta)} =
 {\text{Var}^*_\theta [\overline{f}]}.
\end{eqnarray*}
This, together with (\ref{2eqn:gen_cr_ineq}), yields the following result.
\vspace{0.2cm}

\begin{corollary}
    \label{2cor:gen_cr_inq_for_m_alpha}
    Consider an $\mathbb{M}^{(\alpha)}$-family as in Theorem \ref{2thm:best_est_m_alpha}. Let $\tau^*(\theta) = \mathbb{E}^*_\theta[\bar{f}]$. Let $\widehat{\theta}$ be any unbiased estimator of $\tau^*(\theta)$ with respect to $p_\theta^*$. Then 
\begin{equation}
\label{2eqn:gen_cr_ineq_with_bound}
   \text{Var}_\theta^*[\widehat{\theta}]\geq   \frac{{\Big[\frac{\partial}{\partial\theta}\tau^*(\theta)\Big]^2}}{I^*_{\alpha,n}(\theta)}.
\end{equation}
Equality holds in \eqref{2eqn:gen_cr_ineq_with_bound} if and only if $\widehat{\theta}=\bar{f}$ almost surely.
\end{corollary}
\vspace{0.2cm}

\begin{remark}
\begin{itemize}
    \item [(i)] Observe that
\begin{eqnarray*}
    { \text{Cov}_\theta^*[s^*(X_1^n,\theta),p_\theta^*(X_1^n)^{\alpha -1}s^*(X_1^n,\theta)]^2}
     &\leq &
    \text{Var}_\theta^*[p_\theta^*(X_1^n)^{\alpha -1}s^*(X_1^n,\theta)] \text{Var}_\theta^*[s^*(X_1^n,\theta)],
\end{eqnarray*}
by Cauchy-Schwarz inequality. This implies that
\begin{equation*}
  \dfrac{1}{I_n^*(\theta)}\leq \dfrac{1}{I^*_{\alpha,n}(\theta)}.  
\end{equation*}
That means the bound given by Corollary \ref{2cor:gen_cr_inq_for_m_alpha} is sharper than the usual CRLB for any unbiased estimator of $\tau^*(\theta)$. Observe that Jones et al. estimators are, in general, more robust than MLE but less efficient. The generalized bound in Corollary \ref{2cor:gen_cr_inq_for_m_alpha} might help in investigating the efficacies of such robust estimators.
\vspace{0.2cm}

\item[(ii)] In view of (\ref{2eq:p_star_coincies_p_jones_et_al}) and Remark \ref{2rem:remark_on_m_alpha_family}(i), we see that (\ref{2eqn:gen_cr_ineq_with_bound}) in Corollary \ref{2cor:gen_cr_inq_for_m_alpha} is essentially a generalization of the usual Cram\'er-Rao inequality for an exponential family.
\end{itemize}

\end{remark}
\vspace{0.2cm}

\begin{example}
	Consider the Student distributions with location parameter $\mu=0$, scale parameter $\sigma^2$ (unknown) and a fixed degrees of freedom parameter $\nu\in(2,\infty):$
	\begin{equation}
	\label{student_distr_nu_gmama_form}
	    p_\sigma(x) = N_{\nu,\sigma}\bigg[1+\frac{x^2}{\nu\sigma^2}\bigg]^{-\frac{\nu+1}{2}},
	\end{equation}
	where $N_{\nu,\sigma}$ is the normalizing constant. Let $\alpha = 1-\frac{2}{\nu+1}$, $b_\alpha=1/\nu$ and $N_{\alpha,\sigma} = N_{\nu,\sigma}$. Then (\ref{student_distr_nu_gmama_form}) can be re-written as
a one parameter regular $\mathbb{M}^{(\alpha)}$-family with $h\equiv 1, w(\sigma) = b_\alpha/\sigma^2>0$ and $f(x) = x^2$ \cite{GayenK21J}. By Theorem \ref{2thm:best_est_m_alpha},
	$\overline{x^2}$ is the best estimator for $\mathbb{E}_\sigma^*[\overline{X^2}]$. We shall now find $\mathbb{E}_\sigma^*[\overline{X^2}]$ explicitly.
	Using (\ref{rao_blackwell_distri_jones_et}), we have
	\begin{eqnarray*}
	p^*_\sigma(x_1^n) 
	= Z_{\alpha,\sigma}\Big[1+(b_\alpha/\sigma^2)\overline{x^2}\Big]^{\frac{1}{\alpha -1}},
	\end{eqnarray*}
	where $Z_{\alpha,\sigma}^{-1}\coloneqq \int [1+(b_\alpha/\sigma^2)\overline{y^2}]^{\frac{1}{\alpha -1}} dy_1^n =\sigma^n\cdot H_n{(\alpha)}$, provided $\alpha>1-\frac{2}{n}$. Here,
\begin{itemize}
    \item [] if $n=2k+1$, $k\geq 0$,
    \begin{equation*}
        H_n(\alpha) = \frac{(\pi b_\alpha)^{\frac{1}{2}}(n\pi/b_\alpha)^{\frac{n}{2}}\Gamma\big(\frac{1+\alpha}{2(1-\alpha)}\big)A_n(\alpha)^{k}}{\Gamma(n/2)\Gamma(\frac{1}{1-\alpha})},\quad\text{and}
    \end{equation*}
    \item [] if $n=2k$, $k\geq 1$, 
    \begin{equation*}
        H_n(\alpha) = \frac{(n\pi/b_\alpha)^{\frac{n}{2}}(1-\alpha)A_n(\alpha)^{k-1}}{\alpha\Gamma(n/2)}
    \end{equation*}
\end{itemize}
	with
	$A_n(\alpha)\coloneqq \frac{(1-\alpha)(n-2)}{2-n(1-\alpha)}$. Let us now denote
	$B_n(\alpha)= \frac{n(1-\alpha)}{2\alpha-n(1-\alpha)}$. Then $\mathbb{E}_\sigma^*[\overline{X^2}]$ equals to $\frac{A_n(\alpha)^{2-k}B_n(\alpha)^{k-1}}{b_\alpha}\sigma^2$ when $n=2k$, and $\frac{A_n(\alpha)^{1-k}B_n(\alpha)^{k}}{b_\alpha}\sigma^2$ when $n=2k+1$, provided $\alpha>\frac{n}{n+2}$. Hence the generalized CRLB is given by,
\begin{itemize}
    \item [] when $n=2k+1$,
    \begin{equation*}
        \hspace{-0.5cm} [\{{A_n(\alpha)^{2-k}(C_n(\alpha)^k -A_n(\alpha)^{-k}B_n(\alpha)^{2k}})\}\big/{b_\alpha^2}]\sigma^4,
    \end{equation*}
    \item [] and when $n=2k$,
    \begin{equation*}
        \hspace{-0.3cm} [\{{A_n(\alpha)^{3-k}(C_n(\alpha)^{k-1} -A_n(\alpha)^{1-k}B_n(\alpha)^{(2k-1)}})\}\big/{b_\alpha^2}]\sigma^4
    \end{equation*}
\end{itemize}
where $C_n(\alpha) \coloneqq \frac{(1-\alpha)(n+2)}{2\alpha+(n+2)(\alpha -1)}$ and $\alpha>\frac{n+2}{n+4}$.
\end{example}
\vspace{0.2cm}

\begin{remark}
Observe that Theorem \ref{2thm:best_est_m_alpha} can not be applied, for example, to Student distributions with unknown location parameter (scale parameter known), since they do not form a regular $\mathbb{M}^{(\alpha)}$-family \cite[Rem. 2(c)]{GayenK21J}. However, in the following, we address this problem by deriving a Rao-Blackwell type theorem for Basu et al. likelihood function (\ref{basu_likelihood}). Note that, for this model, Basu et al. estimator and Jones et al. estimator are the same \cite{GayenK18ISIT, GayenK21J}.
\end{remark}

\subsection{Analogous results for Basu et al. likelihood function}
\vspace{0.2cm}

Consider the Basu et al. likelihood function $L^{(\alpha)}_B$ as in (\ref{basu_likelihood}).
The associated deformed probability distribution is given by
\begin{eqnarray}
\label{2eqn:gen_joint_distribution_basu}
    p_\theta^*(x_1^n) &\coloneqq& \frac{\exp[L^{(\alpha)}_B(x_1^n;\theta)]}{\int \exp[L^{(\alpha)}_B(y_1^n;\theta)] d{y_1^n}}\nonumber\\
    &=& \frac{\exp[\frac{\alpha}{n(\alpha -1)}\sum_{i=1}^n p_\theta(x_i)^{\alpha -1}]}{\int \exp[\frac{\alpha}{n(\alpha -1)}\sum_{i=1}^n p_\theta(y_i)^{\alpha -1}]dy_1^n}.~~~~~
\end{eqnarray}
Recall that, $\mathbb{B}^{(\alpha)}$-families are the ones that have a fixed number of sufficient statistics with respect to Basu et al. likelihood function (Theorem \ref{2thm:necessary_condition_sufficiency_basu}).
Consider a $\mathbb{B}^{(\alpha)}$-family as in (\ref{2eqn:form_of_general_b_alpha_family}).
  Then, using (\ref{2eqn:gen_joint_distribution_basu}), we have
\begin{eqnarray}
\label{2eqn:gen_joint_distribution_b_alpha}
p_\theta^*(x_1^n) 
= \exp \Big[\frac{\alpha}{\alpha -1} \overline{h} + \widetilde{F}(\theta) +\frac{\alpha}{\alpha -1}w(\theta) \overline{f}\Big],
\end{eqnarray}
where $\widetilde{F}(\theta) = -\log \int \exp \big[\frac{\alpha}{\alpha -1} (\overline{h}  +w(\theta) \overline{f})\big] dx_1^n$, is the normalizing constant.
Observe that $p_\theta^*$ forms a one parameter exponential family. Then we can derive the following result concerning the best estimator for the parameter, analogous to Theorem \ref{2thm:best_est_m_alpha} and Corollary \ref{2cor:gen_cr_inq_for_m_alpha}.
\vspace{0.2cm}

\begin{theorem}
	\label{2thm:best_est_b_alpha}
Consider a regular $\mathbb{B}^{(\alpha)}$-family as in (\ref{2eqn:form_of_general_b_alpha_family}). Let $\mathbb{E}^*_\theta[\bar{f}] = \tau(\theta)$. Then $\bar{f}$ is the best estimator of its expected value with respect to $p_\theta^*$. Furthermore, if $\widehat{\theta}$ is any estimator  such that $\mathbb{E}^*_\theta[\widehat{\theta}] = \tau(\theta)$, then
	\begin{equation}
\label{2eqn:gen_crlb_b_alpha}
    \text{Var}_\theta^*[\widehat{\theta}] \geq \frac{[\tau'(\theta)]^2}{I_n^*(\theta)}.
\end{equation}
\end{theorem}
\vspace{0.2cm}

\begin{IEEEproof}
    See Sec-D.3 of Appendix.
\end{IEEEproof}
\vspace{0.2cm}

\begin{remark}
Notice that the bound in (\ref{2eqn:gen_crlb_b_alpha}) is given by the usual Fisher information (of course, applied to the model $p_\theta^*$) unlike the one in Corollary \ref{2cor:gen_cr_inq_for_m_alpha} where the bound is given by the asymptotic variance of Jones et al. estimators.
\end{remark}
\vspace{0.2cm}

We now apply the above result to find best estimator for the mean parameter of a Student distribution.
\vspace{0.2cm}

\begin{example}
    Consider the Student distributions with degrees of freedom parameter $\nu\in(2,\infty)$, location parameter $\mu$  and scale parameter $\sigma^2=1:$
\begin{equation}
\label{student_distr_nu}
    p_\mu(x) = N_\nu\Big[1+\frac{1}{\nu}(x-\mu)^2\Big]^{-\frac{\nu+1}{2}},
\end{equation}
where $N_\nu$ is the normalizing constant. Letting $\alpha = 1-\frac{2}{\nu+1}$, $b_\alpha=1/\nu$ and $N_\alpha= N_\nu$, (\ref{student_distr_nu}) forms a regular $\mathbb{B}^{(\alpha)}$-family as in (\ref{2eqn:form_of_general_b_alpha_family}) with
$h(x) = N_\alpha^{\alpha -1}b_\alpha x^2$, $F(\mu) = (1+b_\alpha \mu^2) N_\alpha^{\alpha -1}$, 
$w(\mu) = -2\mu N_\alpha^{\alpha -1}b_\alpha$ and $f(x) = x$.
This implies that $\overline{x}$ is a sufficient statistic for $\mu$ with respect to the Basu et al. likelihood function \cite[Th. 4]{GayenK21ISIT} and $\overline{x}$ is the best estimator for $\mathbb{E}^*_\mu[\overline{X}]$ (Theorem \ref{2thm:best_est_b_alpha}). 

It can also be shown, in this case, that
\begin{eqnarray*}
    {p_\mu^*(x_1^n)} = \exp\Big[\frac{\alpha}{\alpha -1} N_\alpha^{\alpha -1} b_\alpha \overline{x^2} + F_\alpha(\mu) - \frac{2\alpha}{\alpha -1}N_\alpha^{\alpha -1} b_\alpha \mu \overline{x}\Big],
\end{eqnarray*}
where $$F_\alpha(\mu) \coloneqq \frac{\alpha}{\alpha -1}N_\alpha^{\alpha -1}b_\alpha \mu^2 - \frac{n}{2}\log \Big[\pi -  \frac{\alpha}{n(\alpha -1)}N_\alpha^{\alpha -1}b_\alpha\Big],$$
\begin{eqnarray*}
    \mathbb{E}^*_\mu[\overline{X}] &=& -\Big(\frac{\alpha-1}{\alpha}\Big)\frac{\frac{2\alpha}{(1-\alpha)}N_\alpha^{\alpha -1}b_\alpha \mu}{w'(\mu)} = \mu,\quad\text{and}\\
\text{Var}_\mu^*[\overline{X}] &=& \frac{[\tau'(\mu)]^2}{I_n^*(\mu)} = 
\frac{(1-\alpha)}{2\alpha b_\alpha N_\alpha^{\alpha -1}}.
\end{eqnarray*}
\end{example}

\section{Summary and concluding remarks}
\label{2sec:summary}
Fisher introduced the concept of sufficiency in connection with parameter estimation problems in an attempt to estimate the unknown parameters of a population without knowing the entire sample \cite{Fisher22J}. Fisher's notion of sufficiency is based on the usual log-likelihood function and the underlying estimation is ML estimation. 
However, there are situations (for example, robust inference) where we need to look beyond ML estimation. Hence the classical notion of sufficiency is inadequate in such situations. In this paper, we presented a notion of sufficiency based on a general likelihood function.
In particular, we considered two families of likelihood functions, namely Basu et al. and Jones et al. $\alpha$-likelihood functions. These (for $\alpha>1$) arise in distance-based robust inference methods and are generalizations of the usual log-likelihood function. Also, these likelihood functions are special as they do not require smoothing of empirical measure (of the observed sample) \cite{JonesHHB01J}. We characterized the family of probability distributions that have a fixed number of sufficient statistics (independent of sample size) with respect to these likelihood functions. This family generalizes the usual exponential family and has a power law. Student distributions fall in this family. We also extended the notion of minimal sufficiency and found minimal sufficient statistics for the power-law families. This could potentially help when the data is modeled using one of these power-law distributions. However, we would like to emphasize that the established results have the limitation (when it comes to practical applications) that, the tuning parameter $\alpha$ in Basu et al. or Jones et al. likelihood function should be related to the degrees of freedom parameter $\nu$ of the Student distribution by $\alpha={(\nu-1)}/{(\nu+1)}$. Basu et al. and Jones et al. estimations for exponential families, in particular, for normal distributions have been studied, for example, in \cite{BasuHHJ98J, JonesHHB01J, BroniatowskiV12J}. However, unlike these works, the present work focused on models that have closed-form solutions for these estimations.

Another important application of sufficient statistics is that, given an estimator, one can find another estimator that is uniformly better than the given one. This result is known as the Rao-Blackwell theorem. We generalized Rao-Blackwell theorem when one has sufficient statistics with respect to Basu et al. or Jones et al. likelihood function. This helps us find the ``best'' estimator for the parameters of a power-law family and hence for Student distributions. Consequently, we derived a lower bound for the variance of estimators of a deformed form of power-law family, $p_\theta^*$. In case of Jones et al. likelihood function, the bound turns out to be its asymptotic variance (with respect to $p_\theta^*$). The derivation of asymptotic variance of Jones et al. estimator can be found in \cite{JonesHHB01J}. This bound is sharper than the usual CRLB (See Corollary \ref{2cor:gen_cr_inq_for_m_alpha}). Recall that the usual CRLB equals the asymptotic variance of MLE.

We now make some general remarks on the deformed probability distribution $p_\theta^*$. Recall that, it is the factorization theorem (Halmos and Savage \cite{HalmosS49J} in the classical case and Proposition \ref{2prop:gen_factorization_thm} in the general case) that makes it clear that the notion of sufficiency is with respect to a likelihood function. However, the generalized notion of sufficiency is defined using the deformed probability distribution $p_\theta^*$. The equivalence of analogous factorization theorem and Fisher's definition of sufficiency shows the uniqueness of the usage of $p_\theta^*$ in the generalized notion of sufficiency. This fact was further evidenced in the case of Jones et al. likelihood function, where the lower bound of variance of estimators for power-law families equals the asymptotic variance of the Jones et al. estimator for $p_\theta^*$. However, for Basu et al. likelihood function, the bound is given by the usual CRLB, though computed with respect to $p_\theta^*$, does not equal the asymptotic variance of the Basu et al. estimator (Theorem \ref{2thm:best_est_b_alpha}). This is surprising as one would expect the bound to equal the Basu et al. asymptotic variance given by \cite[Th.~2]{BasuHHJ98J} (as in the case of Jones et al. likelihood function). This will be the focus of one of our forthcoming works. In any case, such an anomalous behavior of Basu et al. likelihood function (from Jones et al. likelihood function) was also noticed in the study of the geometry of the associated divergence functions in connection with projection theorems for families of probability distributions determined by linear constraints. While the geometry of the Basu et al. divergence is similar to that of the Kullback-Leibler divergence, the same for Jones et al. divergence was different \cite{KumarS15J2}.


%

\appendix
	\section*{A. Proof of results of Section 2}
	\vspace{0.4cm}
	
\begin{itemize}
    \item [A.1.]

\begin{IEEEproof}[Proof of Proposition \ref{2prop:gen_factorization_thm}]
	\vspace{0.2cm}
	
	Let $T$ be a sufficient statistic for $\theta$
	with respect to $L_G$. Let  $x_1^n$ and $y_1^n$ be two samples such that
	$T(x_1^n) = T(y_1^n)$. 
	Then from generalized Koopman's definition of sufficiency, we have
	$[L_G(x_1^n;\theta) - L_G(y_1^n;\theta)]$ is independent of
	$\theta$. Let us define a relation `$\sim$' on the set of all sample points of length 
	$n$ by $x_1^n\sim y_1^n$ if and only if 
	$T(x_1^n) = T(y_1^n)$.
	Then `$\sim$' is an equivalence relation. Let us denote the equivalence classes by
	$\mathbb{S}^t$ for each $t\in \mathscr{T} = \{r:r=T(x_1^n)~\text{for some}~x_1^n\in\mathbb{S}^n\}$. For each 
	equivalence class $\mathbb{S}^t$, designate an element 
	${x}_{1,t}^n\in \mathbb{S}^t$. 
	
	Let $x_1^n$ be a sample point such that
	$T(x_1^n)= t^*\in \mathscr{T}$. Then $x_1^n\in \mathbb{S}^{t^*}$
	and	$T(x_1^n)= T({x}_{1,t^*}^n)$. This implies that $[L_G(x_1^n;\theta) - 
	L_G({x}_{1,t^*}^n;\theta)]$ is
	independent of $\theta$ by hypothesis. Let $v(x_1^n) := 
	[L_G(x_1^n;\theta) - L_G({x}_{1,t^*}^n;\theta)]$.
	Then
	\begin{eqnarray*}
	    L_G(x_1^n;\theta) 
	= L_G({x}_{1,t^*}^n;\theta) + v(x_1^n)
        = u(\theta, t^*) + v(x_1^n),
	\end{eqnarray*}
	where $u(\theta, t):= L_G({x}_{1,t}^n;\theta)$ for each $t\in \mathscr{T}$.

	Conversely, let us assume that there exist two functions
	$u$ and $v$ such that
	\begin{equation*}
	L_G(x_1^n;\theta) = u (\theta, T(x_1^n)) + 
	v(x_1^n)
	\end{equation*}
	for all sample points $x_1^n$ and for all $\theta\in \Theta$. Thus for any two
	sample points $x_1^n$ and $y_1^n$ with
	$T(x_1^n) = T(y_1^n)$, we have
	\begin{align*}
	{L_G(x_1^n;\theta) - L_G(y_1^n;\theta)}& = u (\theta, T(x_1^n)) + v(x_1^n) 
	- u (\theta, T(y_1^n)) - v(y_1^n)\\
	&= v(x_1^n) - v(y_1^n).
	\end{align*}
	This implies that $[L_G(x_1^n;\theta) - L_G(y_1^n;\theta)]$ is independent of $\theta$.
	By generalized Koopman's definition of sufficiency,
	$T$ is a sufficient statistic for $\theta$ with respect to $L_G$. This completes the proof.
\end{IEEEproof}	
\vspace{0.2cm}

\item[A.2.]
\begin{IEEEproof}[Proof of Proposition \ref{equi_generalized_fisher_generalized_facto}]
	\vspace{0.2cm}
	
    Let $T$ be sufficient for $\theta$ with respect to $L_G$. Then, from Proposition 2, there exist functions $u$ and $v$ such that
    \begin{equation}
    \label{fact_gen}
   L_G(x_1^n;\theta) =u(T(x_1^n),\theta) + v(x_1^n)
\end{equation}
for all $x_1^n\in\mathbb{S}^n$.
We need to only show $p_{\theta{_{X_1^n|T}}}^*(x_1^n|t)$ is independent of $\theta$ when $t=T(x_1^n)$. Observe that if $t\neq T(x_1^n)$, $p_{\theta{_{X_1^n|T}}}^*(x_1^n|t)=0$ which is independent of $\theta$. Using (\ref{fact_gen}), we have
\begin{eqnarray*}
{p_{\theta{_{X_1^n|T}}}^*(x_1^n|t)}
&=& \dfrac{p_\theta^*(x_1^n)}{\int\limits_{A_{T(x_1^n)}} p^*_\theta(y_1^n) d{y_1^n}} \\
&=& \dfrac{\exp[L_G(x_1^n;\theta)]}{\int\limits_{A_{T(x_1^n)}} \exp[L_G(y_1^n;\theta)]dy_1^n}\nonumber\\
&=& \dfrac{\exp[u(T(x_1^n),\theta)]\exp[v(x_1^n)]}{\int\limits_{A_{T(x_1^n)}} \exp[u(T(x_1^n),\theta)]\exp[v(y_1^n)] dy_1^n}\nonumber\\
&=& \dfrac{\exp[v(x_1^n)]}{\int\limits_{A_{T(x_1^n)}} \exp[v(y_1^n)] dy_1^n},
\end{eqnarray*}
which is independent of $\theta$, where the first equality is by the definition of the generalized conditional distribution.

Conversely, let $p_{\theta{_{X_1^n|T}}}^*(x_1^n|t)$ be independent of $\theta$. Let
$p_{\theta{_{X_1^n|T}}}^*(x_1^n|t) \coloneqq v_1(x_1^n)$.
Then
\begin{eqnarray*}
\dfrac{\exp[L_G(x_1^n;\theta)]}{\int\limits_{A_{T(x_1^n)}} \exp[L_G(y_1^n;\theta)]dy_1^n} =v_1(x_1^n).
\end{eqnarray*}
This implies 
\begin{eqnarray}
\label{2eqn:like_fact_expre}
L_G(x_1^n;\theta) = u(T(x_1^n),\theta)+v(x_1^n),
\end{eqnarray}
where $u(T(x_1^n),\theta) = \log \int\limits_{A_{T(x_1^n)}} \exp[L_G(y_1^n;\theta)]dy_1^n$
and $v(x_1^n) = \log v_1(x_1^n)$. Since the above is true for any $x_1^n\in \mathbb{S}^n$, we have the result.
\end{IEEEproof}
\vspace{0.4cm}

\item[A.3.]
\begin{IEEEproof}[Proof of Proposition \ref{1prop:condition_for_minimal_sufficient_statistic}]
	\vspace{0.2cm}
	
	The condition that if $T(x_1^n) = T(y_1^n)$ then  
	$[L_G(x_1^n;\theta) - L_G(y_1^n;\theta)]$ is independent of $\theta$, for all $\theta \in \Theta$ implies that $T$ is a sufficient statistic for $\theta$ with respect to $L_G$, by Definition 1. Thus the only thing
	we need to prove is that $T$ is further minimal. Let $\widetilde{T}$
	be a sufficient statistic such that 
	$\widetilde{T}({x_1^n}) = \widetilde{T}(y_1^n)$.
	Then by Definition 1,
	$[L_G(x_1^n;\theta) - L_G(y_1^n;\theta)]$
	is independent of $\theta$. Hence by hypothesis, we have $T(x_1^n) = T(y_1^n)$. Therefore $T$ is minimal, by definition.
\end{IEEEproof}

\end{itemize}

\section*{B. Proof of results of Section 3}
\vspace{0.4cm}

\begin{itemize}
    \item [B.1.]
\begin{IEEEproof}[Proof of Theorem \ref{2thm:suffiicent_statistics_M_alpha}]
	\vspace{0.2cm}

		Let $\Pi = \mathbb{M}^{(\alpha)}$ as in Definition 6. Let $x_1^n$ and 
		$y_1^n$ be two i.i.d. samples from some $p_\theta\in \mathbb{M}^{(\alpha)}$ such that $T(x_1^n) = T(y_1^n)$.
		Assume that, for $i=1,\dots,s$, 
		\begin{eqnarray*}
		     {\overline{f_i(x_1^n)}/\overline{h(x_1^n)} = \psi _i (T(x_1^n)), 
		~\text{and}~} 
  \overline{f_i(y_1^n)} /\overline{h(y_1^n)} = \psi _i (T(y_1^n)).  
		\end{eqnarray*}
	Then
		\begin{eqnarray}
		\label{2eqn:difference_jones_likelihood}
			\lefteqn{L_J^{(\alpha)}(x_1^n;\theta) - L_J^{(\alpha)}(y_1^n;\theta)}\nonumber\\
			& & \hspace*{-0.5cm} = \frac{1}{\alpha -1}\log\Big[\frac{1}{n}\sum\limits_{i=1}^n p_\theta(x_i)^{\alpha -1}\Big] - \frac{1}{\alpha -1}\log\Big[\frac{1}{n}\sum\limits_{i=1}^n p_\theta(y_i)^{\alpha -1}\Big]\nonumber\\
			& &\hspace*{-0.5cm} = \frac{1}{\alpha -1} \log \frac{Z(\theta)^{\alpha -1}\big\lbrace
			\overline{h(x_1^n)} + w(\theta)^T \overline{f(x_1^n)}\big\rbrace}{Z(\theta)^{\alpha -1}\big\lbrace
			\overline{h(y_1^n)} + w(\theta)^T \overline{f(y_1^n)}\big\rbrace}\nonumber\\
			& & \hspace*{-0.5cm}= \frac{1}{\alpha -1} \log 
			\frac{\overline{h(x_1^n)}}
			{\overline{h(y_1^n)} } +
			\frac{1}{\alpha -1} \log \frac{
			1 + w(\theta)^T \big\lbrace\overline{f(x_1^n)}\big/ \overline{h(x_1^n)}\big\rbrace}{
			1 + w(\theta)^T \big\lbrace \overline{f(y_1^n)}\big/ \overline{h(y_1^n)}\big\rbrace}\\
			& & \hspace*{-0.5cm}= \frac{1}{\alpha -1} \log \big[
			{\overline{h(x_1^n)}}\big/
			{\overline{h(y_1^n)}}\big]\nonumber,
		\end{eqnarray}
		which is independent of $\theta$. The last equality holds since $T(x_1^n) = T(y_1^n)$ and thus 
		$\overline{f_i(x_1^n)}/\overline{h(x_1^n)} = \overline{f_i(y_1^n)}/\overline{h(y_1^n)}$.
		Hence
		$T$ is a sufficient statistic for $\theta$ with respect to the Jones et al. likelihood function $L_J^{(\alpha)}$.
\end{IEEEproof}
\vspace{0.4cm}

\item[B.2.]
\begin{IEEEproof}[Proof of Theorem \ref{2thm:minimal_sufficient_statistics_M_alpha}]
	\vspace{0.2cm}

We only give the proof for $\mathbb{M}^{(\alpha)}$-family. The result for $\mathbb{B}^{(\alpha)}$-family can be proved similarly.
\vspace{0.2cm}

	Let us consider an $\mathbb{M}^{(\alpha)}$-family as in Definition 6. Assume that the family is regular.
	We will only show that the
	statistic 
	\begin{eqnarray*}
	  \overline{f(x_1^n)}\big/\overline{h(x_1^n)}
    =
    \big(\overline{f_1(x_1^n)}\big/\overline{h(x_1^n)},\dots,\overline{f_k(x_1^n)}\big/\overline{h(x_1^n)}\big)
	\end{eqnarray*}
	is minimal. Let us consider two sample points
	$x_1^n$ and $y_1^n$ such that
	$[L_J^{(\alpha)}(x_1^n;\theta) - L_J^{(\alpha)}(y_1^n;\theta)]$ is independent of 
	$\theta$. Our aim is to prove that
	$\overline{f(x_1^n)}\big/\overline{h(x_1^n)}=\overline{f(y_1^n)}\big/\overline{h(y_1^n)}$.
	
	Consider $[L_J^{(\alpha)}(x_1^n;\theta) - L_J^{(\alpha)}(y_1^n;\theta)]$. Using (\ref{2eqn:difference_jones_likelihood}), we have
	\begin{align*}
	{L_J^{(\alpha)}(x_1^n;\theta) - L_J^{(\alpha)}(y_1^n;\theta)
	= \frac{1}{\alpha -1} \log 
			\frac{\overline{h(x_1^n)}}
			{\overline{h(y_1^n)} }}  +\frac{1}{\alpha -1} \log \frac{
			1 + w(\theta)^T \big\lbrace\overline{f(x_1^n)}\big/ \overline{h(x_1^n)}\big\rbrace}{
			1 + w(\theta)^T \big\lbrace \overline{f(y_1^n)}\big/ \overline{h(y_1^n)}\big\rbrace}.
	\end{align*}
	Since the above is independent of $\theta$, we must have 
	\begin{equation*}
	    \log \frac{
			1 + w(\theta)^T \big\lbrace\overline{f(x_1^n)}\big/ \overline{h(x_1^n)}\big\rbrace}{
			1 + w(\theta)^T \big\lbrace \overline{f(y_1^n)}\big/ \overline{h(y_1^n)}\big\rbrace}
			= C,
	\end{equation*}
	where $C$ is a constant independent of $\theta$ (may depend on $x_1^n$ and $y_1^n$). This implies
	\begin{eqnarray*}
	    {(1-\exp(C)) + w(\theta)^T \big[ \overline{f(x_1^n)}\big/ \overline{h(x_1^n)}} - \exp(C)  \big\lbrace \overline{f(y_1^n)}\big/ \overline{h(y_1^n)}\big\rbrace\big] =0.
	\end{eqnarray*}
	
	Since the given family is regular, the functions 1, $w_1,\dots,w_k$ are
	linearly independent. Hence we have $\exp(C) = 1$ and therefore
	\begin{equation*}
	    \big[ \overline{f(x_1^n)}\big/ \overline{h(x_1^n)} -\overline{f(y_1^n)}\big/ \overline{h(y_1^n)} \big]=0.
	\end{equation*}
	That is,
	\begin{equation*}
	   \overline{f(x_1^n)}\big/ \overline{h(x_1^n)} = \overline{f(y_1^n)}\big/ \overline{h(y_1^n)}. 
	\end{equation*}
	Thus $\bar{f}/\bar{h}$ is a minimal sufficient statistic for  a regular $\mathbb{M}^{(\alpha)}$-family with respect to $L_J^{(\alpha)}$.
\end{IEEEproof}
	
\end{itemize}

\section*{C. Proof of results of Section 4}	
\vspace{0.4cm}

\begin{itemize}
    \item[C.1.]
\begin{IEEEproof}[Proof of $\text{Cov}_\theta(\psi(\bar{f}),\bar{f})=0$ for regular exponential family]
\vspace{0.2cm}

Recall that $p_\theta(x) = Z(\theta) [h(x)+w(\theta)f(x)]$ for $x\in\mathbb{S}$ and $\mathbb{E}_\theta[\psi(\bar{f})] = 0$. That is,
\begin{equation*}
\int \psi(\bar{f}) p_\theta(x_1^n) dx_1^n = 0.
\end{equation*}
Taking derivative with respect to $\theta$ on both sides, we get
\begin{equation*}
\int \psi(\bar{f}) \frac{\partial}{\partial\theta} \bigg\{\prod\limits_{i=1}^n\exp[h(x_i) +w(\theta)f(x_i)]\bigg\} dx_1^n =0.
\end{equation*}
This implies
\begin{equation*}
\int \psi(\bar{f}) p_\theta(x_1^n) w'(\theta) n\bar{f} dx_1^n = 0.
\end{equation*}
Since $w'(\theta)\neq 0$ for all $\theta\in\Theta$, we have $\mathbb{E}_\theta[\psi(\bar{f})\bar{f}] = 0.$
Therefore 
\begin{eqnarray*}
   \text{Cov}_\theta[\psi(\bar{f}),\bar{f}]
    = \mathbb{E}_\theta[\psi(\bar{f})\bar{f}] -\mathbb{E}_\theta[\psi(\bar{f})] \mathbb{E}_\theta[\bar{f}] = 0.
\end{eqnarray*}
\end{IEEEproof}

\item[C.2.]
\begin{IEEEproof}[Proof of Proposition \ref{thm1}]
	\vspace{0.2cm}

	Since $T$ is a sufficient statistic with respect to the Jones et al. likelihood function, by Proposition 3, the conditional distribution of  $x_1^n$ given $T$ with respect to $p_\theta^*$ is independent of $\theta$. Hence the result follows as $\widehat{\theta}$ is a function of $x_1^n$ only.
\end{IEEEproof}
\vspace{0.4cm}

\item[C.3.]
\begin{IEEEproof}[Proof of Proposition \ref{thm2}]
	\vspace{0.2cm}

For every partition set $A_{T(y_1^n)}$ of $\mathbb{S}^n$, the space of all possible $n$-length samples, we have
\begin{eqnarray}
\label{expectation_1}
{\mathbb{E}_\theta^*\Big[\phi^*(T(X_1^n))\cdot \mathbf{1}_{A_{T(y_1^n)}}\Big] }
&\coloneqq& \int\limits_{ A_{T(y_1^n)}} p_\theta^*(x_1^n) \phi^*(T(x_1^n)) dx_1^n \nonumber\\
&= &\int\limits_{ A_{T(y_1^n)}} p_\theta^*(x_1^n) \phi^*(T(y_1^n))dx_1^n,
\end{eqnarray}
since $x_1^n\in A_{T(y_1^n)}$ implies $T(x_1^n) = T(y_1^n)$.
\vspace{0.2cm}

Observe that
\begin{eqnarray*}
{\phi^*(T(x_1^n))}
&=& \int \widehat{\theta}(y_1^n) p_\theta^*(y_1^n|T(x_1^n))dy_1^n\nonumber\\
&=& \int\limits_{A_{T(x_1^n)}} \widehat{\theta}(y_1^n) p_\theta^*(y_1^n|T(x_1^n))dy_1^n \nonumber\\
&& + \int\limits_{A^c_{T(x_1^n)}} \widehat{\theta}(y_1^n) p_\theta^*(y_1^n|T(x_1^n))dy_1^n \nonumber\\
&=& \int\limits_{A_{T(x_1^n)}} \widehat{\theta}(y_1^n) p_\theta^*(y_1^n|T(x_1^n))dy_1^n +0 \nonumber\\
&=& \int\limits_{A_{T(x_1^n)}}  \widehat{\theta}(y_1^n) \cdot \Bigg[\dfrac{p_\theta^*(y_1^n)}{q_\theta^*(T(x_1^n))}\Bigg]dy_1^n \nonumber\\
&=& {\int\limits_{A_{T(x_1^n)}} \widehat{\theta}(y_1^n) p_\theta^*(y_1^n)dy_1^n}\Big/{\int\limits_{A_{T(x_1^n)}} p_\theta^*(y_1^n)dy_1^n}.
\end{eqnarray*}

Thus (\ref{expectation_1}) implies
\begin{eqnarray}
\label{expectation_2}
{\mathbb{E}_\theta^*\Big[\phi^*(T(X_1^n))\cdot \mathbf{1}_{A_{T(y_1^n)}}\Big]}
&=& \int\limits_{ A_{T(y_1^n)}} p_\theta^*(x_1^n) \phi^*(T(y_1^n))dx_1^n\nonumber\\
&=&  \phi^*(T(y_1^n)) \int\limits_{ A_{T(y_1^n)}} p_\theta^*(x_1^n)dx_1^n\nonumber\\
&=& \frac{\int\limits_{ A_{T(y_1^n)}} \widehat{\theta}(x_1^n) p_\theta^*(x_1^n) dx_1^n}{\int\limits_{ A_{T(y_1^n)}} p_\theta^*(x_1^n)dx_1^n}
\Big[\int\limits_{ A_{T(y_1^n)}} p_\theta^*(x_1^n)dx_1^n\Big] \nonumber\\
&=& \int\limits_{ A_{T(y_1^n)}} \widehat{\theta}(x_1^n)p_\theta^*(x_1^n) dx_1^n\nonumber\\
&=& \mathbb{E}_\theta^*\Big[\widehat{\theta}(X_1^n))\cdot \mathbf{1}_{A_{T(y_1^n)}}\Big].
\end{eqnarray}
Observe that (\ref{expectation_2}) holds for every partition set $A_{T(y_1^n)}$ of $\mathbb{S}^n$.
This implies both the estimators $\widehat{\theta}$ and  $\phi^*(T)$ have the same expected value with respect to $p_\theta^*$. That is,
\begin{equation*}
  \mathbb{E}_\theta^*[\phi^*(T(X_1^n))] =   \mathbb{E}_\theta^*[\widehat{\theta}(X_1^n))]=\tau^*(\theta).
\end{equation*}
\end{IEEEproof}
\vspace{0.2cm}

\item[C.4.]
\begin{IEEEproof}[Proof of Proposition \ref{thm3}]
	\vspace{0.2cm}

    $\mathcal{A}\subseteq \mathcal{B}$ is immediate from
	Propositions 14 and 15.
	Conversely, let $\psi(T)\in\mathcal{B}$. Since $\mathbb{E}^*_{\theta}[\psi(T)] = \tau^*(\theta)$, we have $\mathbb{E}^*_{\theta}[\psi^*(T)|T]\in \mathcal{A}$. But $\mathbb{E}^*_{\theta}[\psi^*(T)|T]= \psi^*(T)$. This implies that $\psi^*(T)\in\mathcal{A}$. Hence $\mathcal{B}\subseteq\mathcal{A}$. This completes the proof.
\end{IEEEproof}
\vspace{0.4cm}

\item[C.5.]
\begin{IEEEproof} [Proof of Lemma \ref{thm4}]
	\vspace{0.2cm}

We have
	\begin{eqnarray*}
		{\text{Var}^*_{\theta}[\widehat{\theta}]}
        &=& \mathbb{E}^*_{\theta}[(\widehat{\theta} - \tau^*(\theta))^2] \\
		&=& \mathbb{E}^*_{\theta}[(\widehat{\theta} - \psi^*(T) +\psi^*(T) -\tau^*(\theta))^2]\\
		&=& \mathbb{E}^*_{\theta}[(\widehat{\theta} - \psi^*(T))^2] +
		\mathbb{E}^*_{\theta}[(\psi^*(T)-\tau^*(\theta))^2] \\
        &&\hspace{1cm} + 2\mathbb{E}^*_{\theta}[(\widehat{\theta} - \psi^*(T))(\psi^*(T) - \tau^*(\theta))]  \\
		&=&  \mathbb{E}^*_{\theta}[(\widehat{\theta} - \psi^*(T))^2] + 
		\text{Var}^*_{\theta}[\psi^*(T)] \\
        &&\hspace{1cm} + 2\mathbb{E}^*_{\theta}[(\widehat{\theta} - \psi^*(T))(\psi^*(T) - \tau^*(\theta))].
	\end{eqnarray*}
Now \begin{eqnarray*}
	{\mathbb{E}^*_{\theta}[(\widehat{\theta} - \psi^*(T))(\psi^*(T) - \tau^*(\theta))] }
	&=& \mathbb{E}^*_{\theta}[(\widehat{\theta} - \psi^*(T))\psi^*(T)]  -\tau^*(\theta)~ \mathbb{E}^*_{\theta}[(\widehat{\theta} - \psi^*(T))]\\
	&\stackrel{(a)}{=}&  \mathbb{E}^*_{\theta}[(\widehat{\theta} - \psi^*(T))\psi^*(T)]\\
	&=& \mathbb{E}^*_{\theta}[\widehat{\theta}~ \psi^*(T)] -\mathbb{E}^*_{\theta}[(\psi^*(T))^2]\\
	&\stackrel{(b)}{=}& \mathbb{E}^*_{\theta}[\mathbb{E}^*_{\theta}[\widehat{\theta}~ \psi^*(T)|T]] -\mathbb{E}^*_{\theta}[(\psi^*(T))^2]\\
	&=&  \mathbb{E}^*_{\theta}[\psi^*(T)~\mathbb{E}^*_{\theta}[\widehat{\theta} |T]] -\mathbb{E}^*_{\theta}[(\psi^*(T))^2]\\
	&=& \mathbb{E}^*_{\theta}[\psi^*(T)\big\{\mathbb{E}^*_{\theta}[\widehat{\theta} |T] - \psi^*(T)\big\}]\\
	&=& \mathbb{E}^*_{\theta}[\psi^*(T)\{\phi^*(T) - \psi^*(T)\}],
\end{eqnarray*}
where (a) follows from the fact that both $\widehat{\theta}$ and $\psi^*(T)$ are unbiased estimators of $\tau^*(\theta)$ and (b) follows from the law of total expectation.
Thus we have 
\begin{eqnarray*}
    \text{Var}^*_{\theta}[\widehat{\theta}] &=& \mathbb{E}^*_{\theta}[(\widehat{\theta} - \psi^*(T))^2] + \text{Var}^*_{\theta}[\psi^*(T)] + 2 \mathbb{E}^*_{\theta}[\psi^*(T)\{\phi^*(T) -\psi^*(T)\}].
\end{eqnarray*}
\end{IEEEproof}
\vspace{0.4cm}

\item[C.6.]
\begin{IEEEproof}[Proof of Theorem \ref{2thm:gen_rao_blackwell}]
	\vspace{0.2cm}

	From Proposition 15, we have $\mathbb{E}^*_{\theta}[\phi^*(T)] = \mathbb{E}^*_{\theta}[\widehat{\theta}] = \tau^*(\theta)$. Now in Lemma 17, if we take $\psi^*(T)= \phi^*(T)$, we have
	\begin{eqnarray*}
		\text{Var}^*_{\theta}[\widehat{\theta}] &=& \mathbb{E}^*_{\theta}[(\widehat{\theta} - \phi^*(T))^2] + \text{Var}^*_{\theta}[\phi^*(T)]\\
  &\geq&  \text{Var}^*_{\theta}[\phi^*(T)].
	\end{eqnarray*}
Equality holds if and only if $\mathbb{E}^*_{\theta}[(\widehat{\theta} - \phi^*(T))^2] =0$, that is, if and only if $\widehat{\theta} = \phi^*(T) = \mathbb{E}^*_\theta[\widehat{\theta}|T]$ almost surely. This happens if and only if $\widehat{\theta}$ is a function of $T$ only. This completes the proof.
\end{IEEEproof} 
\vspace{0.4cm}

\item[C.7.]
\begin{IEEEproof}[Proof of Proposition \ref{unique_best_est_wrt_p_star}]
\vspace{0.2cm}

   Let $\phi_1(T)$ and $\phi_2(T)$ be two unbiased estimators for $\tau^*(\theta)$ with respect to $p_\theta^*$. Let $\sigma^2_*$ be the minimum variance that any unbiased estimator for $\tau^*(\theta)$ can achieve. Suppose $\text{Var}_\theta^*[\phi_1(T)] = \text{Var}_\theta^*[\phi_2(T)]= \sigma^2_*$. We show that $\phi_1(T) = \phi_2(T)$ almost surely.

    Consider $\xi(T):= \dfrac{\phi_1(T) + \phi_2(T)}{2}$. Then $\mathbb{E}_\theta^*[\xi(T)] = \tau^*(\theta)$ and 
    \begin{eqnarray*}
    \text{Var}^*_\theta[\xi(T)] \leq \sigma_*^2,
    \end{eqnarray*}
    since $\text{Cov}_\theta^*[\phi_1(T), \phi_2(T)] \leq \sigma_*^2$, by Cauchy-Schwarz inequality. As $\sigma^2_*$ is the minimum variance that one unbiased estimator of $\tau^*(\theta)$ can achieve, we have $\text{Var}^*_\theta[\xi(T)] \geq \sigma_*^2$. Therefore $\text{Var}^*_\theta[\xi(T)] = \sigma_*^2$. This implies
    \begin{equation*}
        \text{Cov}_\theta^*[\phi_1(T), \phi_2(T)] = \sigma_*^2.
    \end{equation*}
    That means $\phi_1(T)$ and $\phi_2(T)$ are positively correlated. Thus there exist constants, say $\alpha>0$ and $\beta$, such that
    \begin{equation*}
        \phi_1(T) =\alpha \phi_2(T) + \beta\quad\text{almost surely.}
    \end{equation*}
    Taking expectation and variance with respect to $p_\theta^*$ on both sides, respectively, we get
    \begin{equation*}
        \tau^*(\theta) = \alpha \tau^*(\theta) + \beta\quad \text{and}\quad \sigma_*^2 = \alpha^2 \sigma_*^2.
    \end{equation*}
    These together imply $\alpha =1$ and $\beta=0$. That is, $\phi_1(T)=\phi_2(T)$ almost surely.
\end{IEEEproof}
\end{itemize}

\section*{D. Proof of results of Section 5}
\vspace{0.4cm}

\begin{itemize}
    \item [D.1.]
   \begin{IEEEproof}[Proof of $\text{Cov}_\theta^*(\psi^*(\bar{f}),\bar{f})\geq 0$ in a regular $\mathbb{M}^{(\alpha)}$-family]
\vspace{0.2cm}

Since $\mathbb{E}_\theta^*[\psi^*(\bar{f})] = 0$,
\begin{equation*}
\int \psi^*(\bar{f}) N(\theta) [1+w(\theta)\overline{f}]^{\frac{1}{\alpha -1}} dx_1^n =0.
\end{equation*}
That is,
\begin{equation*}
\int \psi^*(\bar{f}) [1+w(\theta)\overline{f}]^{\frac{1}{\alpha -1}} dx_1^n =0.
\end{equation*}
Differentiating both sides with respect to $\theta$,
\begin{equation*}
\int \psi^*(\bar{f}) \frac{1}{\alpha -1} [1+w(\theta)\overline{f}]^{\frac{2-\alpha}{\alpha -1}} w'(\theta) \bar{f}dx_1^n =0.
\end{equation*}
This implies
\begin{equation*}
\int\frac{\psi^*(\bar{f})\bar{f}}{[1+w(\theta)\overline{f}]} p_\theta^*(x_1^n) dx_1^n = 0.
\end{equation*}
That is,
\begin{equation}
\label{2eqn:unbiased_est_zero}
\mathbb{E}^*_\theta\Bigg[\frac{\psi^*(\bar{f})\bar{f}}{1+w(\theta)\overline{f}}\Bigg] = 0.
\end{equation}
We now claim that $\mathbb{E}_\theta^*[\psi^*(\bar{f})\bar{f}] \geq 0$.
Let
\begin{equation}
\label{2eqn:def-of_k_theta}
 k(\theta)\coloneqq \mathbb{E}_\theta^*[\psi^*(\bar{f})\bar{f}].
\end{equation}
 Suppose that $k(\theta)\neq 0$ for all $\theta\in\Theta$. From (\ref{2eqn:def-of_k_theta}),
\begin{equation*}
\int \psi^*(\bar{f}) \bar{f} N(\theta)[1+w(\theta)\bar{f}]^{\frac{1}{\alpha -1}} dx_1^n = k(\theta).
\end{equation*}
That is,
\begin{equation*}
\int \psi^*(\bar{f}) \bar{f} [1+w(\theta)\bar{f}]^{\frac{1}{\alpha -1}} dx_1^n = \frac{k(\theta)}{N(\theta)}.
\end{equation*}
Differentiating both sides with respect to $\theta$, we get
\begin{eqnarray*}
   \int \psi^*(\bar{f}) \bar{f}\frac{1}{\alpha -1} [1+w(\theta)\bar{f}]^{\frac{2-\alpha}{\alpha -1}} w'(\theta)\bar{f}dx_1^n
   &=& \frac{\partial}{\partial\theta}\Big[\frac{k(\theta)}{N(\theta)}\Big].
\end{eqnarray*}
Hence
\begin{eqnarray*}
    {\int \psi^*(\bar{f})\frac{{\bar{f}}^2}{1+w(\theta)\bar{f}} p^*_\theta(x_1^n)dx_1^n }
&=& (\alpha -1) \frac{N(\theta)}{w'(\theta)} \frac{\partial}{\partial\theta}\Big[\frac{k(\theta)}{N(\theta)}\Big].
\end{eqnarray*}
Hence
\begin{eqnarray*}
    {\int \psi^*(\bar{f}) \bar{f} \Big[\frac{1+w(\theta)\bar{f} -1}{1+w(\theta)\bar{f}}\Big] p_\theta^*(x_1^n) dx_1^n}
    &=& (\alpha -1) \frac{w(\theta)N(\theta)}{w'(\theta)} \frac{\partial}{\partial\theta}\Big[\frac{k(\theta)}{N(\theta)}\Big].
\end{eqnarray*}
This implies
\begin{eqnarray*}
    {\mathbb{E}_\theta^*[\psi^*(\bar{f})\bar{f}] - \mathbb{E}^*_\theta\Big[\frac{\psi^*(\bar{f})\bar{f}}{1+w(\theta)\overline{f}}\Big]}
&=& (\alpha -1) \frac{w(\theta)N(\theta)}{w'(\theta)} \frac{\partial}{\partial\theta}\Big[\frac{k(\theta)}{N(\theta)}\Big].
\end{eqnarray*}
Using (\ref{2eqn:unbiased_est_zero}), we get
\begin{equation*}
k(\theta) = (\alpha -1) \frac{w(\theta)N(\theta)}{w'(\theta)} \frac{\partial}{\partial\theta}\Big[\frac{k(\theta)}{N(\theta)}\Big].
\end{equation*}
Since $k(\theta)\neq 0$, we have
\begin{equation*}
(\alpha-1)\frac{\partial}{\partial\theta}\Big[\frac{k(\theta)}{N(\theta)}\Big]\Big/ \Big[\frac{k(\theta)}{N(\theta)}\Big] = \frac{w'(\theta)}{w(\theta)}.
\end{equation*}
That is,
\begin{equation*}
(\alpha -1)\frac{\partial}{\partial\theta}\log \Big[\frac{k(\theta)}{N(\theta)}\Big] = \frac{\partial}{\partial\theta} \log w(\theta).
\end{equation*}
This reduces to
\begin{equation*}
\frac{\partial}{\partial\theta}\log \Big[\frac{k(\theta)}{N(\theta)}\Big]^{\alpha -1} = \frac{\partial}{\partial\theta} \log w(\theta).
\end{equation*}
Integrating both sides with respect to $\theta$, we get
\begin{equation*}
\log \Big[\frac{k(\theta)}{N(\theta)}\Big]^{\alpha -1} = \log w(\theta) + \log C
\end{equation*}
where $\log C$ is the constant of integration.
This implies
\begin{equation*}
k(\theta)^{\alpha -1} = C\cdot w(\theta)N(\theta)^{\alpha -1}.
\end{equation*}
Observe that $C>0$ and $N(\theta)>0$ for all $\theta\in\Theta$. Hence, if $w(\theta) >0$ for all $\theta\in\Theta$, we get $k(\theta)> 0$. This, together with the case $k(\theta)=0$, implies $k(\theta)\geq 0$.
That means, in any case, $\mathbb{E}^*_\theta[\psi^*(\bar{f})\bar{f}] \geq 0$. Therefore
\begin{eqnarray*}
    \text{Cov}_\theta^*[\psi^*(\bar{f}),\bar{f}]
    &=& \mathbb{E}^*_\theta[\psi^*(\bar{f})\bar{f}] - \mathbb{E}^*_\theta[\psi^*(\bar{f})]
\mathbb{E}^*_\theta[\bar{f}]\\
&=& \mathbb{E}^*_\theta[\psi^*(\bar{f})\bar{f}] -0\\
&=&\mathbb{E}^*_\theta[\psi^*(\bar{f})\bar{f}]\\
&\geq& 0.
\end{eqnarray*}
\end{IEEEproof}
\vspace{0.2cm}

\item[D.2.]
\begin{IEEEproof}[Derivation of generalized Cram\'er-Rao bound for $\mathbb{M}^{(\alpha)}$]
\vspace{0.2cm}

\noindent
Consider 
\begin{eqnarray*}
p_\theta^*(x_1^n)
    &=& N(\theta) [1+w(\theta)\overline{f}]^{\frac{1}{\alpha -1}}\\
    &=& [N(\theta)^{\alpha -1} + N(\theta)^{\alpha -1}w(\theta) \overline{f}]^{\frac{1}{\alpha -1}}.
\end{eqnarray*}
\noindent
For convenience, let $F(\theta) := N(\theta)^{\alpha -1}$ and $\widetilde{w}(\theta) := N(\theta)^{\alpha -1}w(\theta)$. Then
\begin{equation*}
    p_\theta^*(x_1^n)
    = [F(\theta) + \widetilde{w}(\theta) \overline{f}]^{\frac{1}{\alpha -1}},
\end{equation*}
and thus
\begin{eqnarray}
    \label{2eq:derivative_of_p_star}
   \frac{\partial}{\partial\theta}p_\theta^*(x_1^n)
    = \frac{1}{\alpha -1} p_\theta^*(x_1^n)^{2-\alpha} [F'(\theta) + \widetilde{w}'(\theta) \overline{f}].~~
\end{eqnarray}
Hence, from the definition of $s^*$,
\begin{eqnarray*}
\lefteqn{\text{Cov}^*_\theta[s^*(X_1^n,\theta), p_\theta^*(X_1^n)^{\alpha -1}s^*(X_1^n,\theta)]}\nonumber\\
&& \hspace{-0.6cm}= \int s^*(x_1^n,\theta)^2 p^*_\theta (x_1^n)^\alpha dx_1^n\nonumber\\
&& \hspace{-0.6cm} = \int s^*(x_1^n,\theta)\Big[\frac{1}{p^*_\theta(x_1^n)}\frac{\partial}{\partial \theta}p^*_\theta(x_1^n)\Big] p^*_\theta(x_1^n)^\alpha dx_1^n\nonumber\\
&&\hspace{-0.6cm}  = \tfrac{1}{\alpha -1} \int p^*_\theta(x_1^n) s^*(x_1^n,\theta) [F'(\theta) +\widetilde{w}'(\theta) \bar{f}] dx_1^n\nonumber\\
&& \hspace{-0.6cm} = \frac{1}{\alpha-1} \int p^*_\theta(x_1^n) s^*(x_1^n,\theta) F'(\theta) dx_1^n + 
\frac{1}{\alpha-1} \int p^*_\theta(x_1^n) s^*(x_1^n,\theta) \widetilde{w}'(\theta) \bar{f}dx_1^n\nonumber\\
&&\hspace{-0.6cm}  =  \frac{1}{\alpha-1} \widetilde{w}'(\theta)\int p^*_\theta(x_1^n) \frac{1}{p^*_\theta(x_1^n)} \frac{\partial}{\partial\theta}p^*_\theta(x_1^n) \bar{f} dx_1^n\nonumber\\
&&\hspace{-0.6cm} = \frac{1}{\alpha-1} \widetilde{w}'(\theta)\frac{\partial}{\partial\theta} \int  p^*_\theta(x_1^n) \bar{f} dx_1^n\nonumber\\
&&\hspace{-0.4cm}  = \frac{1}{\alpha-1} \widetilde{w}'(\theta) \frac{\partial}{\partial\theta}\tau^*(\theta).
\end{eqnarray*}

Using (\ref{2eq:derivative_of_p_star}), we get
\begin{align*}
 {\int s^*(x_1^n,\theta) p_\theta^*(x_1^n)^\alpha dx_1^n}
& = \int \frac{\partial}{\partial\theta} \log p_\theta^* (x_1^n) p_\theta^*(x_1^n)^\alpha dx_1^n\nonumber\\
&= \frac{1}{\alpha -1} \int p_\theta^*(x_1^n) [F'(\theta) +\widetilde{w}'(\theta)\bar{f}] dx_1^n
\end{align*}
and
\begin{align*}
{\int s^*(x_1^n,\theta)^2 p_\theta^*(x_1^n)^{2\alpha -1} dx_1^n}
&= \int \Big[\frac{\partial}{\partial\theta}\log p_\theta^*(x_1^n)\Big]^2 p_\theta^*(x_1^n)^{2\alpha -1} dx_1^n \nonumber\\
&=\Big(\frac{1}{\alpha -1}\Big)^2 \int p_\theta^*(x_1^n) [F'(\theta) +\widetilde{w}'(\theta)\bar{f}]^2 dx_1^n.
\end{align*}
Thus
\begin{align*}
{\text{Var}_\theta^*[p_\theta^*(X_1^n)^{\alpha -1}s^*(X_1^n,\theta)]}
& = \int s^*(x_1^n,\theta)^2 p_\theta^*(x_1^n)^{2\alpha -1} dx_1^n - \Big[\int s^*(x_1^n,\theta) p_\theta^*(x_1^n)^\alpha dx_1^n\Big]^2\\
& = \Big(\frac{1}{\alpha -1}\Big)^2 \{\mathbb{E}^*_\theta[(F'(\theta) +\widetilde{w}'(\theta)\bar{f})^2] -(\mathbb{E}^*_\theta [F'(\theta) +\widetilde{w}'(\theta)\bar{f}])^2\}\nonumber\\
&= \Big(\frac{1}{\alpha -1}\Big)^2 \text{Var}^*_\theta [F'(\theta) +\widetilde{w}'(\theta)\bar{f}]\nonumber\\
&= \Big(\frac{1}{\alpha -1}\Big)^2 [\widetilde{w}'(\theta)]^2 \text{Var}^*_\theta [\bar{f}].
\end{align*}
\end{IEEEproof}
\vspace{0.2cm}

\item[D.3.]
\begin{IEEEproof}[Proof of Theorem \ref{2thm:best_est_b_alpha}]
\vspace{0.2cm}

The deformed probability distribution associated with Basu et al. likelihood function is given by
\begin{eqnarray}
\label{deform_dist_basu}
   p_\theta^*(x_1^n)
    = \dfrac{\exp[\frac{\alpha}{n(\alpha -1)}\sum_{i=1}^n p_\theta(x_i)^{\alpha -1}]}{\int \exp[\frac{\alpha}{n(\alpha -1)}\sum_{i=1}^n p_\theta(y_i)^{\alpha -1}]dy_1^n}.~~
\end{eqnarray}
$T:=T(x_1^n)$ is a sufficient statistic for $\theta$ with respect to Basu et al. likelihood function. Let $\widehat{\theta}:=\widehat{\theta} (x_1^n)$ be an unbiased estimator of $\tau(\theta)$ with respect to $p_\theta^*$.  Let us define
\begin{eqnarray*}
     \phi(T(x_1^n))
     &:=& \mathbb{E}_\theta^*[\widehat{\theta}|T=T(x_1^n)]\\
     &=& \int \widehat{\theta}(y_1^n) p_\theta^*(y_1^n|T(x_1^n))dy_1^n,
\end{eqnarray*}
where the conditional distribution is as defined in (11) of the main article. Then proceeding with similar steps as in Lemma 17 and Theorem 18, we have the following.

\begin{itemize}
    \item [1.] $\mathbb{E}_\theta^* [\phi(T)] = \mathbb{E}_\theta^*[\widehat{\theta}]= \tau(\theta)$,
    
    \item[2.] $\text{Var}_\theta^*[\phi(T)]\leq \text{Var}_\theta^*[\widehat{\theta}]$, and equality holds if and only if $\widehat{\theta}$ is a function of $T$ only.
\end{itemize}

Now for the given $\mathbb{B}^{(\alpha)}$-family, using (\ref{deform_dist_basu}), we have 
\begin{equation*}
p_\theta^*(x_1^n) 
= \exp \Big[\frac{\alpha}{\alpha -1} \overline{h} + \widetilde{F}(\theta) +\frac{\alpha}{\alpha -1}w(\theta) \overline{f}\Big].
\end{equation*}
 Recall that $T = \overline{f}$ is a sufficient statistic for $\theta$ with respect to Basu et al. likelihood function. Using generalized Rao-Blackwell theorem for Basu et al. likelihood function, we get
\begin{eqnarray}
    \label{2eqn:rao_blackwell_b_alpha}
  &&\mathbb{E}_\theta^*[\phi(\overline{f})] = \mathbb{E}_\theta^*[\widehat{\theta}] =\tau(\theta)\quad\text{and}\nonumber\\
    &&\text{Var}_\theta^*[\phi(\bar{f})]\leq \text{Var}_\theta^*[\widehat{\theta}].
\end{eqnarray}
Let us define $\psi(\bar{f})\coloneqq \phi(\bar{f}) -\bar{f}$. Then $\mathbb{E}_\theta^*[\psi(\overline{f})] = 0$, since $\mathbb{E}_\theta^*[\bar{f}] = \tau(\theta)$ by hypothesis.
Also
\begin{eqnarray*}
   {\text{Var}_\theta^*[\phi(\bar{f})]}
    &=& \text{Var}_\theta^*[\psi(\bar{f}) +\bar{f}]\\
    &=& \text{Var}_\theta^*[\psi(\bar{f})]
    + \text{Var}_\theta^*[\bar{f}] + 2~ \text{Cov}^*(\psi(\bar{f}), \bar{f}).
\end{eqnarray*}
Proceeding as in C.1., one can show that $\text{Cov}^*(\psi(\bar{f}), \bar{f}) = 0$. Thus we get
\begin{equation*}
    \text{Var}_\theta^*[\phi(\bar{f})] = \text{Var}_\theta^*[\psi(\bar{f})]
    + \text{Var}_\theta^*[\bar{f}] \geq \text{Var}_\theta^*[\bar{f}].
\end{equation*}
Hence the inequality in (\ref{2eqn:rao_blackwell_b_alpha}) reduces to
\begin{equation}
    \text{Var}_\theta^*[\widehat{\theta}] \geq \text{Var}_\theta^*[\bar{f}],
\end{equation}
that is, $\bar{f}$ is the best estimator of its expected value with respect to $p_\theta^*$. 
\vspace{0.2cm}

We next show that $\text{Var}_\theta^*[\bar{f}] = \dfrac{[\tau'(\theta)]^2}{I_n^*(\theta)}$.
Observe that $p_\theta^*$ forms a one-parameter exponential family and
\begin{eqnarray}
\label{2eqn:fisher_inf_p_star}
I_n^*(\theta)&=& \text{Var}_\theta^*[\frac{\partial}{\partial\theta}\log p_\theta^*(X_1^n)]\nonumber\\
&=& \text{Var}_\theta^*[\widetilde{F}'(\theta) + \frac{\alpha}{\alpha -1}w'(\theta)\bar{f}]\nonumber\\
&=& \Big(\frac{\alpha}{\alpha -1}\Big)^2 (w'(\theta))^2~ \text{Var}_\theta^*[\bar{f}].
\end{eqnarray}
Using the fact $\mathbb{E}^*_\theta[\frac{\partial}{\partial\theta}\log p_\theta^*(X_1^n)] =0$ and $\mathbb{E}_\theta^*[\Bar{f}] = \tau(\theta)$, we get
\begin{equation}
\label{2eqn:exp_value_of_f_bar}
    \tau(\theta) = -\frac{(\alpha -1)}{\alpha}\frac{\widetilde{F}'(\theta)}{w'(\theta)}.
\end{equation}
Further taking derivative on both sides of $\mathbb{E}^*_\theta[\frac{\partial}{\partial\theta}\log p_\theta^*(X_1^n)] =0$, we get 
\begin{equation}
\label{2eqn:derivative_exp_value_of_f_bar}
    \tau'(\theta) = \frac{\alpha}{\alpha -1} w'(\theta) \text{Var}_\theta^*[\bar{f}].
\end{equation}
Thus, using (\ref{2eqn:fisher_inf_p_star}), (\ref{2eqn:exp_value_of_f_bar}) and (\ref{2eqn:derivative_exp_value_of_f_bar}), we get
\begin{equation*}
     \text{Var}_\theta^*[\overline{f}] = \frac{[\tau'(\theta)]^2}{I_n^*(\theta)}.
\end{equation*}
This completes the proof.
\end{IEEEproof}
\end{itemize}




\section*{Acknowledgments}
The authors would like to thank the Editor, Associate Editor, and Referees for their valuable comments that improved the presentation of the paper.

\bibliographystyle{IEEEtranS} 
\bibliography{myjmva}       

\end {document}